\documentclass[11pt]{article}
\usepackage{amsmath}
\usepackage{amscd}
\usepackage{amsthm}
\usepackage{verbatim}
\usepackage{fullpage}
\usepackage{amssymb}
\usepackage{amsfonts}
\usepackage{mathrsfs}
\usepackage{stmaryrd}
\usepackage[usenames]{color}
\usepackage{graphicx}
\usepackage{hyperref}
\usepackage{mathscinet}
\usepackage{palatino}
\usepackage{bbm}
\usepackage{latexsym}
\usepackage{exscale}
\usepackage[normalem]{ulem}

\usepackage[backend=biber,style=alphabetic,doi=false,isbn=false,url=false,minnames=6,maxnames=6,maxalphanames=6]{biblatex}
\renewbibmacro{in:}{%
  \ifentrytype{article}{}{\printtext{\bibstring{in}\intitlepunct}}}
\renewbibmacro*{volume+number+eid}{%
  \printtext{vol.}
  \printfield{volume}%
  \setunit*{\addnbspace}
  \printfield{number}%
  \setunit{\addcomma\space}%
  \printfield{eid}}
\DeclareFieldFormat[article]{number}{\mkbibparens{#1}}

\newcommand{\googlebooks}[1]{(preview at \href{https://books.google.com/books?id=#1}{google books})}

\newcommand{\numdam}[1]{}

\theoremstyle{plain}
\newtheorem{proposition}{Proposition}

\numberwithin{equation}{section}
\newtheorem{example}[proposition]{Example}

\newtheorem{remark}[proposition]{Remark}

\usepackage[all]{xy}
\usepackage{graphicx}
\usepackage{tikz}
\usepackage{tikz-cd}
\usetikzlibrary{calc}
\usetikzlibrary{decorations.markings}
\usetikzlibrary{decorations.pathreplacing}
\usetikzlibrary{arrows,shapes,positioning}
\tikzstyle directed=[postaction={decorate,decoration={markings,
    mark=at position #1 with {\arrow{>}}}}]
\tikzstyle rdirected=[postaction={decorate,decoration={markings,
    mark=at position #1 with {\arrow{<}}}}]
\tikzset{anchorbase/.style={baseline={([yshift=-0.5ex]current bounding box.center)}},
  tinynodes/.style={font=\tiny,text height=0.75ex,text depth=0.15ex},
  smallnodes/.style={font=\scriptsize,text height=0.75ex,text depth=0.15ex},
  >={Latex[length=1mm, width=1.5mm]}
}

\tikzset{
    partial ellipse/.style args={#1:#2:#3}{
        insert path={+ (#1:#3) arc (#1:#2:#3)}
    }
}
\newcommand{\ru}{to [out=0,in=270]}
\newcommand{\rd}{to [out=0,in=90]}
\newcommand{\ur}{to [out=90,in=180]}

\newcommand{\ul}{to [out=90,in=0]}

\newcommand{\dr}{to [out=270,in=180]}
\newcommand{\dl}{to [out=270,in=0]}
\newcommand{\pu}{to [out=90,in=270]}

\def\R{\mathbb{R}}
\def\Z{\mathbb{Z}}

\def\C{\mathbb{C}}

\newcommand{\kb}{{\mathbbm k}}

\newcommand{\glnn}[1]{\mathfrak{gl}_{#1}}
\newcommand{\slnn}[1]{\mathfrak{sl}_{#1}}
\newcommand{\glN}{\mathfrak{gl}_N} 
\newcommand{\slN}{\mathfrak{sl}_N}
\newcommand{\gltw}{\mathfrak{gl}_2} 
\newcommand{\sltw}{\mathfrak{sl}_2}
\newcommand{\Uqsltw}{U_q(\mathfrak{sl}_2)}
\newcommand{\Uqgltw}{U_q(\mathfrak{gl}_2)}
\newcommand{\UqslN}{U_q(\mathfrak{sl}_N)} 
\newcommand{\UqglN}{U_q(\mathfrak{gl}_N)}

\DeclareMathOperator{\id}{id}

\newcommand{\CC}[1]{\left\llbracket #1 \right\rrbracket}

\newcommand{\Sym}{\mathrm{Sym}}
\newcommand{\KhR}{\mathrm{KhR}_N}
\newcommand{\KhRnn}[1]{\mathrm{KhR}_{#1}}
\newcommand{\KhRS}{\mathrm{KhR}_{\Sigma}}

\newcommand{\Kh}{\KhR}

\newcommand{\skeinzero}{\mathcal{S}^N_{0}} 

\definecolor{kwcolor}{rgb}{.05, .5, .3}
\definecolor{pwcolor}{rgb}{.15, .5, .5}
\definecolor{dred}{rgb}{.7, 0, 0}

\title{A note on general linear link homology} 
\author{Paul Wedrich\footnote{Universität Hamburg, Fachbereich Mathematik, Bundesstraße 55, 20146 Hamburg, Germany
email: {\tt paul.wedrich@uni-hamburg.de}}}

\begin{document}
\maketitle
\begin{abstract} This expository note outlines why it is sometimes useful to consider the
bigraded type A link homology theories as associated with the Lie algebras
$\glN$ instead of $\slN$.
\end{abstract}

Khovanov's categorification of the Jones polynomial~\cite{MR1740682} is the
paradigmatic example of a link homology theory. It associates chain complexes of
graded abelian groups to link diagrams and chain maps to elementary transitions
between link diagrams, such as Reidemeister moves, creation and annihilation of
unknot components, and saddle moves. Khovanov homology is a categorification in
the sense that the graded Euler characteristic of the chain complex for a link
diagram computes the associated Jones polynomial, up to a certain normalisation.
\smallskip

The Jones polynomial can be identified with the Reshetikhin--Turaev
\cite{MR1036112} link invariant associated to the quantum group $\Uqsltw$, with
the vector representation coloring all link components. The analogous link
invariants for the family of type A quantum groups $\UqslN$ for $N\geq 2$ are
also known as specialisations of the 2-variable HOMFLYPT
invariant~\cite{MR776477,MR945888}.

Categorifications of the quantum link invariants for $\UqslN$ have first been
constructed using matrix factorizations by Khovanov--Rozansky~\cite{MR2391017}
and subsequently in many other contexts. The resulting
invariants are often referred to as $\slN$ \emph{link homology theories}. 
\smallskip

The purpose of this note is to outline why it is sometimes useful to consider
the finite-rank type $A$ link homologies as associated with the Lie algebra
$\glN$ instead of $\slN$. For the sake of an overview, we list three aspects here
that are discussed in the following.

\begin{enumerate}
  \item Constructions of finite-rank type A link homologies that use an integral
  quantum grading are forced to categorify \emph{a} $\UqglN$ skein relation
  instead of \emph{the} $\UqslN$ skein relation.
  \item Deformations of the underlying commutative Frobenius algebra lead to
  spectral sequences that relate the finite-rank type A link homologies to
  filtered link homologies. The latter can be decomposed
  along branching rules for $\glN$-representations. For example, Khovanov homology for
  $\gltw$ deforms into Lee homology for $\glnn{1}\oplus \glnn{1}$.
  \item Making Khovanov homology functorial under link cobordisms requires
  breaking the symmetry of the complex associated to a crossing of two (unoriented)
  strands, e.g. by using foams which encode the non-triviality of
  the determinant representation of $\gltw$. This also fixes sign-issues in base point actions and enables gradings by intergral homology in related skein module constructions.
\end{enumerate}

\subsection*{Conventions}
  For the sake of concreteness, we will work with the combinatorially constructed
  finite-rank type $A$ link homologies that are based on the closed foam
  evaluation formula of Robert--Wagner~\cite{1702.04140}. The construction and the
  proof of functoriality under link cobordisms appears in Ehrig--Tubbenhauer--Wedrich
  \cite{MR3877770}. Several other constructions of type A link homologies are
  proven to produce isomorphic chain complexes when evaluated on
  links\footnote{Such isomorphisms are rarely canonical and not much is known
  about functoriality under link cobordisms in these other constructions.} and
  much of our discussion here carries over to these invariants. 
We denote the $\glN$ Khovanov--Rozansky homology by $\KhR$. The invariant of the unknot is
\begin{equation*}
 \KhR(\bigcirc) \cong q^{1-N}\C[X]/\langle X^N\rangle \cong H^{*-N+1}(\C P^{N-1},\C).
\end{equation*}

Here and throughout, $X$ is a variable of (\emph{quantum}) degree $2$ and
$q^{1-N}$ indicates a downward degree shift of magnitude $N-1$. The ground ring
will typically be $\C$ or any field of characteristic zero. 

\subsection*{Acknowledgements}
The author would like to thank Thang Le for a discussion concerning skein
relations for $\glN$ and $\slN$, and Matt Hogancamp for suggesting to compare
the framing dependence of symmetrically- and antisymmetrically colored link
homologies. 

\subsection*{Funding}
The author acknowledges support from the Deutsche
Forschungsgemeinschaft (DFG, German Research Foundation) under Germany's
Excellence Strategy - EXC 2121 ``Quantum Universe'' - 390833306 and the
Collaborative Research Center - SFB 1624 ``Higher structures, moduli spaces and
integrability'' - 506632645.

\section{Integrality of the quantum grading and the type A skein relations}
\label{sec:skeinrel}
We consider the full subcategory of complex representations of $\glN$, which contains the fundamental representations
and is closed under tensor products, finite direct sums and taking direct
summands. As a consequence of Schur--Weyl duality, this can be identified with
the free symmetric monoidal Karoubian $\C$-linear category generated by one
object $V$ \emph{of rank at most} $N$. The rank condition simply states the
vanishing $\wedge^{N+i}(V)=0$ for $i>1$ of sufficiently high anti-symmetric
Schur functors (defined using the symmetric monoidal structure and idempotent
completeness). This extends to a description of the full subcategory generated
by fundamental representations and their duals, if one requries the determinant
$\wedge^{N}(V)$ to be tensor invertible. 

Analogous characterisations are available for monoidal Karoubian subcategories
of $\UqglN$ over $\C(q)$, although with a non-symmetric braiding (at generic
values of $q$). In particular, these categories are still generated under tensor
product (and duality) by the vector representation $V=\C(q)^N$. The braiding 
\begin{equation*}
  \label{eq:braiding}
  \sigma^{SW}_{V\otimes V}\colon  V\otimes V \to V \otimes V   
\end{equation*}
is determined (up to possible inversion and unit scaling, see
\eqref{eq:rescbraiding}) by the action of standard generator $T_1$ of the Hecke
$H_2$ algebra of type $A_1$ under quantum Schur--Weyl duality. When using the
convention\footnote{Here we have $T_1^2= (q-q^{-1})T_1 + 1$. Another common
convention replaces $q=v^{-1}$ in this relation, see e.g. \cite{MR4220642}.}
\[
H_2=\Z[q,q^{-1},T_1]/\langle (T_1-q)(T_1+q^{-1}) \rangle  
\]
the operator of multiplying by $T_1$, and thus $\sigma^{SW}_{V\otimes V}$, have eigenvalues
$q$ and $-q^{-1}$. We denote the corresponding eigenspaces of
$\sigma^{SW}_{V\otimes V}$ by $\wedge^2_q(V)$ and $\Sym^2_q(V)$ respectively. The projection
\begin{equation}
  \label{eq:proj}
  V\otimes V \to \wedge^2_q(V) \hookrightarrow V\otimes V  
\end{equation}
is given by the action of the following element of $H_2\otimes_{\Z[q,q^{-1}]} \C(q)$:
\[ \frac{T_1 - (-q)^{-1}}{q- (-q)^{-1}} = \frac{T_1+q^{-1}}{q+q^{-1}} \] In the numerator we
see the Kazhdan--Lusztig basis element~\cite{MR560412} $B_{1}:= T_1+q^{-1}=T_1^{-1}+q$, which is a
quasi-idempotent with eigenvalue $q+q^{-1}=[2]$. To summarise, the linear
relations between $T_1$, $T_1^{-1}$, $B_{1}$, and the unit are given by:

\begin{equation}
  \label{eq:Heckerel}
T_1 = B_{1}-q^{-1}, \qquad T_1^{-1} = B_{1} - q ,\qquad T_1-T_1^{-1}=q-q^{-1} 
\end{equation}

\subsection{Diagrammatics for type A representation categories}
We will encode the Schur--Weyl action of the Kazhdan--Lusztig basis
element $B_{1}$ diagrammatically by:
\begin{equation*}
  \begin{tikzpicture}[anchorbase,scale=.5,tinynodes]
    \draw [very thick] (0,0) \pu (0,0.1) \pu (.5,.6);
    \draw [very thick] (1,0) \pu (1,0.1) \pu (.5,.6);
    \draw [very thick, directed=.7] (.5,.6) \pu (.5,.9);
    \node at (.8,.75) {$2$};
    \draw [very thick,->] (.5,.9) \pu (0,1.3) \pu (0,1.5);
    \draw [very thick,->] (.5,.9) \pu (1,1.3) \pu (1,1.5);
  \end{tikzpicture} 
  \;\;:=\;\;  \sigma^{SW}_{V,V}\;\;+\;\; q^{-1} \id_{V\otimes V}
  \;\;:=\;\; (\sigma^{SW}_{V,V})^{-1} \;\;+\;\;  q \id_{V\otimes V}
\end{equation*}

The label $2$ at the central edge indicates that we are consider a
(quasi-idempotent) endomorphism of $V\otimes V$ that factors through
$\wedge^2_q(V)$, see \eqref{eq:proj}. More generally, one can define higher
exterior powers $\wedge^k_q(V)$ and diagrammatically encode intertwiners:

\[
\begin{tikzpicture}[scale =.5,anchorbase,tinynodes]
	\draw[very thick] (0,0) node[right,xshift=-2pt]{$a$} \pu (.5,1);
	\draw[very thick] (1,0) node[right,xshift=-2pt]{$b$} \pu (.5,);
	\draw[very thick,->] (.5,1) to (.5,1.5) node[right,xshift=-1pt]{$a{+}b$};
\end{tikzpicture}
\colon  \wedge^a_q(V) \otimes \wedge^b_q(V) \to \wedge^{a+b}_q(V)
,\qquad
\begin{tikzpicture}[scale =.5,rotate=180,anchorbase,tinynodes]
	\draw[very thick,<-] (0,0) node[right,xshift=-2pt]{$b$} \pu (.5,1);
	\draw[very thick,<-] (1,0) node[right,xshift=-2pt]{$a$} \pu  (.5,1);
	\draw[very thick] (.5,1) to (.5,1.5) node[right,xshift=-2pt]{$a{+}b$};
\end{tikzpicture}
\colon \wedge^{a+b}_q(V) \to \wedge^a_q(V) \otimes \wedge^b_q(V)
\]
The planar diagrammatic calculus of (tensor) composites of such morphisms has
been introduced by Murakami--Ohtsuki--Yamada~\cite{MR1659228} and a complete set
of relations follows from the work of
Cautis--Kamnitzer--Morrison~\cite{MR3263166}, see
Tubbenhauer--Vaz--Wedrich~\cite{MR3709658}. Acknowleding Kuperberg's pioneering
work \cite{MR1403861}, such diagrams are often called webs. 

The Kazhdan--Lusztig basis element
of $H_n$ corresponding to the longest element $w_0$ in $S_n$ is given by the
$q$-analog of the Young symmetrizer:
\[B_{w_0}:=q^{-l(w_0)} \sum_{w\in S_n} q^{l(w)} T_w\] where $T_w$ is the
standard basis element of the Hecke algebra associated to the permutation $w\in
S_n$ and $l(w)$ denotes the length of a minimal expression of $w$ in terms of
simple transpositions. Diagrammatically, the image of $B_{w_0}$ under
Schur--Weyl duality is represented as the quasi-idempotent endomorphism
\[
  \begin{tikzpicture}[scale =.5,anchorbase,tinynodes]
    \draw[very thick] (0,0)  \pu (1,1);
    \draw[very thick] (.5,0) node[right,yshift=1pt]{$\cdots$} \pu (1,1) node[right,xshift=-1pt]{$n$};
    \draw[very thick] (2,0)  \pu (1,1);
    \draw[very thick,directed=.55] (1,1) to (1,1.5);
    \draw[very thick,->] (1,1.5)  \pu (0,2.5);
    \draw[very thick,->] (1,1.5)  \pu (.5,2.5) node[right,yshift=-1pt]{$\cdots$};
    \draw[very thick,->] (1,1.5)  \pu (2,2.5);
  \end{tikzpicture}
  \colon V^{\otimes n} \to \wedge^n_q(V) \hookrightarrow V^{\otimes n}
\]
with eigenvalue $[n]!=[n][n-1]\cdots[2]$. The simple fact that
$\wedge^n_q(V)\cong 0$ for $n>N$ is encoded diagrammatically by setting
declaring that every web with an edge of label $n>N$ represents the zero
morphism.\medskip

Let $c\in \C(q)$ be a unit and consider the braiding obtained by rescaling our
previously chosen Schur--Weyl braiding on the monoidal Karoubian generator $V$:
\begin{equation}
  \label{eq:rescbraiding}
  \sigma_{V,V} := c\cdot \sigma^{SW}_{V,V}\end{equation} 
   We
represent $\sigma_{V,V}$ by the diagrammatic crossing:
\begin{equation*}
  \label{eq:diagrel}
  \begin{tikzpicture}[anchorbase,scale=.5]
    \draw [very thick, ->] (1,0) \pu (0,1.5);
    \draw [white, line width=.15cm] (0,0) \pu (1,1.5);
    \draw [very thick, ->] (0,0) \pu (1,1.5);
  \end{tikzpicture}
  \;\;:=\;\; c \left( \begin{tikzpicture}[anchorbase,scale=.5]
    \draw [very thick] (0,0) \pu (0,0.1) \pu (.5,.6);
    \draw [very thick] (1,0) \pu (1,0.1) \pu (.5,.6);
    \draw [very thick, directed=.7] (.5,.6) \pu (.5,.9);
    \node at (.8,.75) {\tiny $2$};
    \draw [very thick,->] (.5,.9) \pu (0,1.3) \pu (0,1.5);
    \draw [very thick,->] (.5,.9) \pu (1,1.3) \pu (1,1.5);
  \end{tikzpicture}  \;\; -\;\;  q^{-1}\;
  \begin{tikzpicture}[anchorbase,scale=.5]
    \draw [very thick, ->] (0,0) -- (0,1.5);
    \draw [very thick, ->] (1,0) -- (1,1.5);
  \end{tikzpicture}
  \right) 
\end{equation*}
The relations inherited from \eqref{eq:Heckerel} are:

\begin{gather*}
  \label{eq:diagrel2}
  \begin{tikzpicture}[anchorbase,scale=.5]
    \draw [very thick, ->] (1,0) \pu (0,1.5);
    \draw [white, line width=.15cm] (0,0) \pu (1,1.5);
    \draw [very thick, ->] (0,0) \pu (1,1.5);
  \end{tikzpicture}
  \;\;:=\;\; c \left( \begin{tikzpicture}[anchorbase,scale=.5]
    \draw [very thick] (0,0) \pu (0,0.1) \pu (.5,.6);
    \draw [very thick] (1,0) \pu (1,0.1) \pu (.5,.6);
    \draw [very thick, directed=.7] (.5,.6) \pu (.5,.9);
    \node at (.8,.75) {\tiny $2$};
    \draw [very thick,->] (.5,.9) \pu (0,1.3) \pu (0,1.5);
    \draw [very thick,->] (.5,.9) \pu (1,1.3) \pu (1,1.5);
  \end{tikzpicture}  \;\; -\;\;  q^{-1}\;
  \begin{tikzpicture}[anchorbase,scale=.5]
    \draw [very thick, ->] (0,0) -- (0,1.5);
    \draw [very thick, ->] (1,0) -- (1,1.5);
  \end{tikzpicture}
  \right) 
  ,\qquad
  \begin{tikzpicture}[anchorbase,xscale=-.5,yscale=.5]
    \draw [very thick, ->] (1,0) \pu (0,1.5);
    \draw [white, line width=.15cm] (0,0) \pu (1,1.5);
    \draw [very thick, ->] (0,0) \pu (1,1.5);
  \end{tikzpicture}
  \;\;:=\;\; c^{-1} \left( \begin{tikzpicture}[anchorbase,scale=.5]
    \draw [very thick] (0,0) \pu (0,0.1) \pu (.5,.6);
    \draw [very thick] (1,0) \pu (1,0.1) \pu (.5,.6);
    \draw [very thick, directed=.7] (.5,.6) \pu (.5,.9);
    \node at (.8,.75) {\tiny $2$};
    \draw [very thick,->] (.5,.9) \pu (0,1.3) \pu (0,1.5);
    \draw [very thick,->] (.5,.9) \pu (1,1.3) \pu (1,1.5);
  \end{tikzpicture}  \;\; -\;\;  q\;
  \begin{tikzpicture}[anchorbase,scale=.5]
    \draw [very thick, ->] (0,0) -- (0,1.5);
    \draw [very thick, ->] (1,0) -- (1,1.5);
  \end{tikzpicture}
  \right)  \\
  c^{-1} \;
  \begin{tikzpicture}[anchorbase,scale=.5]
    \draw [very thick, ->] (1,0) \pu (0,1.5);
    \draw [white, line width=.15cm] (0,0) \pu (1,1.5);
    \draw [very thick, ->] (0,0) \pu (1,1.5);
  \end{tikzpicture}
  \;\;-\;\;
  c\;
  \begin{tikzpicture}[anchorbase,xscale=-.5,yscale=.5]
    \draw [very thick, ->] (1,0) \pu (0,1.5);
    \draw [white, line width=.15cm] (0,0) \pu (1,1.5);
    \draw [very thick, ->] (0,0) \pu (1,1.5);
  \end{tikzpicture}
  \;\;=\;\;
  (q-q^{-1}) \; \begin{tikzpicture}[anchorbase,scale=.5]
    \draw [very thick, ->] (0,0) -- (0,1.5);
    \draw [very thick, ->] (1,0) -- (1,1.5);
  \end{tikzpicture}
\end{gather*}

\noindent The (rescaled) braiding of higher exterior power representations $\sigma_{\wedge^a_q(V),\wedge^b_q(V)}$ is then expressed as

\begin{equation}
\label{eq:gencrossing}
\begin{tikzpicture}[anchorbase,scale=1,tinynodes]
  \draw [very thick, ->] (1,0)node[right,xshift=-2pt]{$b$}  \pu (0,1.5)node[right,xshift=-2pt]{$b$};
  \draw [white, line width=.15cm] (0,0) \pu (1,1.5);
  \draw [very thick, ->] (0,0)node[right,xshift=-2pt]{$a$}  \pu (1,1.5)node[right,xshift=-2pt]{$a$};
\end{tikzpicture}
=
c^{ab}\sum_{k\geq 0} (-q)^{k-b}\;
\begin{tikzpicture}[anchorbase,scale=1,tinynodes]
  \draw [very thick,directed=.55, ->] (1,0)node[right,xshift=-2pt]{$b$}  \pu (1,1.5)node[right,xshift=-2pt]{$a$};
  \draw [very thick] (1,.25) \pu (0,.75);
  \draw [very thick] (0,.75)  \pu (1,1.25);
  \draw [very thick,directed=.55, ->] (0,0)node[right,xshift=-2pt]{$a$}  \pu (0,1.5)node[right,xshift=-2pt]{$b$};
  \node at (.5,0.3) {$k$};
\end{tikzpicture} \qquad a\geq b
\end{equation}
and the inverse braiding is obtained by inverting $q\leftrightarrow q^{-1}$ and $c\leftrightarrow c^{-1}$.

\subsection{Restricting to \texorpdfstring{$\UqslN$}{Uq(sl(N))}}
\label{sec:repslN}
Now let us restrict from $\UqglN$ to the subalgebra $\UqslN$. The only major
difference is that the determinant now becomes trivial, i.e. isomorphic to the
tensor unit. Let us fix such an isomorphism $\phi$ and introduce diagrammatic representatives:
\[
  \begin{tikzpicture}[anchorbase,scale=1,tinynodes]
  \draw [very thick,directed=.55] (1,0)node[right,xshift=-2pt]{$N$}  \pu (1,.5);
  \node at (1,.5) {$\bullet$};
\end{tikzpicture}=\phi\colon \wedge_q^N(V) \xrightarrow{\cong} \C(q) 
\qquad  
\begin{tikzpicture}[anchorbase,scale=1,tinynodes]
  \draw [very thick,->] (1,0)  \pu (1,.5) node[right,xshift=-2pt]{$N$};
  \node at (1,0) {$\bullet$};
\end{tikzpicture}=\phi^{-1}\colon \C(q) \xrightarrow{\cong} \wedge_q^N(V)\quad 
\]
As usual we omit drawing unitors or identity morphisms on the tensor unit. Let
us assume that the rescaled braiding from \eqref{eq:rescbraiding} descends to a
braiding for the representation category of $\UqslN$. Then we use the naturality of the braiding and
\eqref{eq:gencrossing} to compute:

\begin{equation}
  \begin{tikzpicture}[anchorbase,scale=1,tinynodes]
    \draw [very thick, ->] (1,0)node[right,xshift=-2pt]{$b$}  \pu (0,1.5)node[right,xshift=-2pt]{$b$};
    \draw [very thick, directed=.55] (0,0)node[right,xshift=-2pt]{$N$} to [out=90,in=225] (.25,.5)node{$\bullet$};
    \draw [very thick, ->] (.75,1)node{$\bullet$} to [out=45,in=270] (1,1.5)node[right,xshift=-2pt]{$N$};
  \end{tikzpicture}
   = 
   \begin{tikzpicture}[anchorbase,scale=1,tinynodes]
    \draw [very thick, ->] (1,0)node[right,xshift=-2pt]{$b$}  \pu (0,1.5)node[right,xshift=-2pt]{$b$};
    \draw [white, line width=.15cm] (0,0) \pu (1,1.5);
    \draw [very thick, directed=.55] (0,0)node[right,xshift=-2pt]{$N$} to [out=90,in=225] (.25,.5)node{$\bullet$};
    \draw [very thick, ->] (.35,.6)node{$\bullet$} to [out=45,in=270] (1,1.5)node[right,xshift=-2pt]{$N$};
  \end{tikzpicture}
  =
  \begin{tikzpicture}[anchorbase,scale=1,tinynodes]
    \draw [very thick, ->] (1,0)node[right,xshift=-2pt]{$b$}  \pu (0,1.5)node[right,xshift=-2pt]{$b$};
    \draw [white, line width=.15cm] (0,0) \pu (1,1.5);
    \draw [very thick, ->] (0,0)node[right,xshift=-2pt]{$N$}  \pu (1,1.5)node[right,xshift=-2pt]{$N$};
  \end{tikzpicture}
  =
  c^{Nb}\sum_{k\geq 0} (-q)^{k-b}\;
  \begin{tikzpicture}[anchorbase,scale=1,tinynodes]
    \draw [very thick,directed=.55, ->] (1,0)node[right,xshift=-2pt]{$b$}  \pu (1,1.5)node[right,xshift=-2pt]{$N$};
    \draw [very thick] (1,.25) \pu (0,.75);
    \draw [very thick] (0,.75)  \pu (1,1.25);
    \draw [very thick,directed=.55, ->] (0,0)node[right,xshift=-2pt]{$N$}  \pu (0,1.5)node[right,xshift=-2pt]{$b$};
    \node at (.5,0.3) {$k$};
  \end{tikzpicture} 
  =
  c^{Nb}(-q)^{-b}\;
  \begin{tikzpicture}[anchorbase,scale=1,tinynodes]
    \draw [very thick,directed=.55, ->] (1,0)node[right,xshift=-2pt]{$b$}  \pu (1,1.5)node[right,xshift=-2pt]{$N$};
    \draw [very thick] (0,.75)  \pu (1,1.25);
    \draw [very thick, ->] (0,0)node[right,xshift=-2pt]{$N$}  \pu (0,1.5)node[right,xshift=-2pt]{$b$};
  \end{tikzpicture} 
  \end{equation}
  However, an analogous relation using the inverse braiding identifies the
  morphism on the left-hand side with the $c^{-Nb}(-q)^{b}$ multiple of the web
  on the right-hand side. An easy computation shows the the relevant morphism is non-zero, and so we deduce the relations
  \[ c^{Nb}(-q)^{-b} = c^{-Nb}(-q)^{b}\;\; 1\leq b\leq N \quad   \leftrightarrow \quad
  c^N=\pm q  \] 
In this way, we see that an $\UqslN$-compatible rescaling of the braiding requires the existence of an $N$th root of $\pm q$. For a choice of $q^{1/N}$, the rescaling by $c=-q^{1/N}$ is diagrammatically generated by

  \begin{equation}
    \label{eq:slskeinrel}
    \begin{tikzpicture}[anchorbase,scale=.5]
      \draw [very thick, ->] (1,0) \pu (0,1.5);
      \draw [white, line width=.15cm] (0,0) \pu (1,1.5);
      \draw [very thick, ->] (0,0) \pu (1,1.5);
    \end{tikzpicture}
    \;\;:=\;\;  -q^{1/N}\;\; 
    \begin{tikzpicture}[anchorbase,scale=.5]
      \draw [very thick] (0,0) \pu (0,0.1) \pu (.5,.6);
      \draw [very thick] (1,0) \pu (1,0.1) \pu (.5,.6);
      \draw [very thick, directed=.7] (.5,.6) \pu (.5,.9);
      \node at (.8,.75) {\tiny $2$};
      \draw [very thick,->] (.5,.9) \pu (0,1.3) \pu (0,1.5);
      \draw [very thick,->] (.5,.9) \pu (1,1.3) \pu (1,1.5);
    \end{tikzpicture} 
    \;\; + \;\; 
    q^{-1+1/N}\;
    \begin{tikzpicture}[anchorbase,scale=.5]
      \draw [very thick, ->] (0,0) -- (0,1.5);
      \draw [very thick, ->] (1,0) -- (1,1.5);
    \end{tikzpicture}
  \end{equation}
  and the specialization of \eqref{eq:gencrossing} at $c=-q^{1/N}$ recovers (the
  inverse) braiding described in \cite[Corollary 6.2.3]{MR3263166}. The sign in
  this rescaling is arbitrary but convenient, since it matches the
  interpretation of $\wedge^2_q(V)$ as a deformation of the exterior power/anti-symmetrization $\wedge^2(V)$ at
  $q=1$.

  \smallskip

  In categorifications of type $A$ link polynomials, the local data associated to a
  crossing is typically a 2-term chain complex of the form\footnote{Often with
  differential directed in the opposite direction. Also, in the original
  construction of Khovanov homology the dumbbell web gets replaced by a cup-cap
  diagram, but this does not change the argument.} 
    \begin{equation}
      \label{eq:crossingcx}
      \CC{
      \begin{tikzpicture}[anchorbase,scale=.45]
        \draw [very thick, ->] (1,0) \pu (0,1.5);
        \draw [white, line width=.15cm] (0,0) \pu (1,1.5);
        \draw [very thick, ->] (0,0) \pu (1,1.5);
      \end{tikzpicture}
      }
      \;\;:=\;\; q^k t^l \left( 0 \to 
      \uwave{\CC{
        \begin{tikzpicture}[anchorbase,scale=.45]
          \draw [very thick] (0,0) \pu (0,0.1) \pu (.5,.6);
          \draw [very thick] (1,0) \pu (1,0.1) \pu (.5,.6);
          \draw [very thick, directed=.7] (.5,.6) \pu (.5,.9);
          \node at (.8,.75) {\tiny $2$};
          \draw [very thick,->] (.5,.9) \pu (0,1.3) \pu (0,1.5);
          \draw [very thick,->] (.5,.9) \pu (1,1.3) \pu (1,1.5);
        \end{tikzpicture} 
        }} \to q^{-1}
        \CC{
      \begin{tikzpicture}[anchorbase,scale=.45]
        \draw [very thick, ->] (0,0) -- (0,1.5);
        \draw [very thick, ->] (1,0) -- (1,1.5);
      \end{tikzpicture}
        } \to 0
        \right)
    \end{equation}  
    over a $\Z$-graded $\kb$-linear category that categorifies the MOY calculus
    for type A representation categories. Here monomials in $q$ and $t$ indicate
    shifts in the integral quantum and cohomological gradings and the
    \uwave{highlighted} term is placed in cohomological degree zero. The
    magnitudes $k,l$ of the global shifts are dependent on conventions. However,
    for lack of fractional $q$-grading shifts, this categorification strategy is
    incompatible with a decategorification to the $\UqslN$-braiding \eqref{eq:slskeinrel}. 

    The convention \eqref{eq:crossingcx} for $k=l=0$ is common when considering Rouquier
    complexes of Soergel bimodules. Rouquier~\cite{0409593} works in this convention when proving that the
    complexes thus associated to braid diagrams only depend on the underlying
    braid up to canonical homotopy equivalence. Moreover,
    Libedinsky--Williamson~\cite{MR3283617} proved that the complexes of
    positive (resp. negative) braid lifts from the symmetric group provide a
    semi-orthogonal decomposition of the dg category of chain complexes of
    Soergel bimodules. For an application, see
    Gorsky--Hogancamp--Wedrich~\cite[Theorem 1.3]{2002.06110}. 
When constructing $\glN$ link homologies using this convention, one obtains
    an invariant of framed link diagrams, for which positive curl in a
    $k$-colored strand can be removed at the expense of a grading shift $q^{-k(N-k+1)} t^{k}$.

    Note that the rescaling of the generating diagrammatic crossing
    \eqref{eq:rescbraiding} forces a rescaling of all crossings as in
    \eqref{eq:gencrossing} to satisfy the naturality constraint of the braiding.
    Naturality is also expected in the categorified setting and proven in a
    prototypical example in \cite{liu2024braided,stroppel2024braidingtypesoergelbimodules}.

\section{Deformations and branching rules}
\label{sec:def}

Let $\Sigma=\{\lambda_1^{N_1},\dots,\lambda_l^{N_l}\}$ denote a multiset of
pairwise distinct complex numbers $\lambda_i$ with multiplicites $N_i$ for
$1\leq i\leq l$ with $\sum_{i=1}^lN_i=N$. Let $P_\Sigma$ denote the monic
polynomial of degree $N$ with root multiset $\Sigma$: 
\[P_\Sigma(X):=(X-\lambda_1)^{N_1}\cdots (X-\lambda_l)^{N_l}\] 

The $\Sigma$-deformed Khovanov--Rozansky homology, denoted $\KhRS$, is a
singly-graded (by cohomological degree) and filtered (instead of the quantum
degree) link homology with unknot invariant:
\begin{equation*}
  \KhRS(\bigcirc) \cong \kb[X]/\langle P_{\Sigma}(X)\rangle 
 \end{equation*}
Motivated by Lee's deformation~\cite{MR2173845} of Khovanov homology, the
deformed analogs $\KhRS$ of Khovanov--Rozansky invariants have been studied by
Gornik~\cite{0402266}, Rasmussen~\cite{MR3447099},
Mackaay--Vaz~\cite{MR2336253}, Wu~\cite{MR3392963},
Rose--Wedrich~\cite{MR3590355}, Robert--Wagner~\cite{1702.04140} and others. \smallskip

The main theorem of \cite{MR3590355} (transported to the construction of
\cite{1702.04140}) gives an interpretation of $\KhRS$ in terms of
(singly-graded) Khovanov--Rozansky homologies of smaller rank. For the unknot invariant this follows
immediately from the Chinese Remainder theorem:
\begin{equation}
  \label{eq:CRT}
  \KhRS(\bigcirc) 
  \cong \kb[X]/\langle P_{\Sigma}(X)\rangle 
  \cong \bigoplus_{i=1}^l \kb[X]/\langle (X-\lambda_i)^{N_i}\rangle 
  \cong \bigoplus_{i=1}^l \KhRnn{N_i}(\bigcirc)
\end{equation}
In fact, such a (singly-graded) isomorphism $\KhRS(K) \cong \bigoplus_{i=1}^l
\KhRnn{N_i}(K)$ exists for every knot $K$, and a suitable generalisation holds
for links, see \cite[Theorem 1.1]{MR3590355}. For every $M\geq 0$ we denote the
vector representation of $\glnn{M}$ by $V_M:=\C^M$ and observe an analogy
between \eqref{eq:CRT} and the branching rule $V_N \cong \bigoplus_{i=1}^l
V_{N_i}$ describing the decomposition of $V_N$ when restricting to the
subalgebra $\bigoplus_{i=1}^l \glnn{N_i}\subset \glN$. 

More generally, if we denote the $\wedge_q^k$-colored version of the knot $K$ by
$K^k$, then \cite[Theorem 1.1]{MR3590355} shows that $\KhRS(K^k)$ decomposes
along the branching rule for $\wedge^k(V_N)$:
\begin{equation}
  \label{eq:coldef}
  \KhRS(K^k) 
  \cong \bigoplus_{\sum k_i=k} \bigotimes_{i=1}^l \KhRnn{N_i}(K^{k_i}), \qquad \wedge^k(V_N) \cong \bigoplus_{\sum k_i=k} \bigotimes_{i=1}^l \wedge^{k_i}(V_{N_i})
\end{equation}

In the following, we illustrate the decomposition behaviour of the deformed Lee
homologies of the Hopf link for $\sltw$ and $\gltw$. These examples show that
the summands of the Lee homology for $\gltw$ have a clear interpretation in
terms of the branching rule. 

\subsection{Lee homology for \texorpdfstring{$\sltw$}{sl(2)} and
\texorpdfstring{$\gltw$}{gl(2)}} 
Lee's deformation of Khovanov homology \cite{MR2173845} arises
by deforming the underlying commutative Frobenius algebra (the unknot invariant)
\[\Kh(\bigcirc) \cong \kb[X]/\langle X^2\rangle \quad \xrightarrow{}\quad  \mathrm{Lee}(\bigcirc) \cong \kb[X]/\langle X^2-1\rangle  \]  
Since $X^2-1=(X-1)(X+1)$, we can describe an analog of Lee homology in the
Khovanov--Rozansky context as $\KhRS$ with $\Sigma=\{+1,-1\}$.
The purpose of this section is to compare $\mathrm{Lee}$ and $\KhRS$ and
illustrate their differences in an example. \smallskip

\begin{example}
For the standard diagram of the positive Hopf link, the cube of resolutions takes the following form\footnote{We use slightly non-standard conventions to match with the $\glN$-behaviour. Analogous arguments apply to all common choices of conventions.}
  \begin{equation}
  \xy
  (-18,0)*{
    \CC{\begin{tikzpicture} [xscale=-.4,yscale=.4,anchorbase]
   \draw[very thick, directed=.6] (1.7,1.3) to [out=315,in=180] (2,1) to [out=0,in=180] (3,2) to [out=0,in=270] (3.5,2.5) to [out=90,in=0] (3,3) to (1,3) to [out=180,in=90] (0.5,2.5) to [out=270,in=180] (1,2) to [out=0,in=135] (1.3,1.7);
   \draw[very thick,directed=.4] (2.7,1.3) to [out=315,in=180] (3,1) to [out=0,in=90] (3.5,0.5) to [out=270,in=0] (3,0) -- (1,0) to [out=180,in=270] (.5,0.5) to [out=90,in=180] (1,1) to [out=0,in=180] (2,2) to [out=0,in=135] (2.3,1.7);
  \end{tikzpicture}} =};
  (80,0)*{
    \begin{tikzpicture} [scale=.4]
      \draw[very thick] (3.5,2.5) to [out=90,in=0] (3,3) to (1,3) to [out=180,in=90] (0.5,2.5);
      \draw[very thick,black] (3.5,0.5) to [out=270,in=0] (3,0) -- (1,0) to [out=180,in=270] (.5,0.5);
       \draw[very thick, black] (0.5,2.5) to [out=270,in=180] (1,2) to (3,2) to [out=0,in=270] (3.5,2.5);
      \draw[very thick,black]  (.5,0.5) to [out=90,in=180] (1,1) to (3,1) to [out=0,in=90] (3.5,0.5);
     \end{tikzpicture}
  };
  (88,0)*{
  \begin{tikzpicture} [scale=.4]
   \draw[opacity=0] (0,0) to (0,-5);
  \end{tikzpicture}
  };
  (20,5.5)*{\begin{tikzpicture} [scale=.5]
   \draw[->] (0,0) -- (1.5,0.5);
    \end{tikzpicture}};
   (20,-5.5)*{\begin{tikzpicture} [scale=.5]
   \draw[->] (0,0) -- (1.5,-0.5);
    \end{tikzpicture}};
   (60,5.5)*{\begin{tikzpicture} [scale=.5]
   \draw[->] (0,0) -- (1.5,-0.5);
    \end{tikzpicture}};
   (60,-5.5)*{\begin{tikzpicture} [scale=.5]
   \draw[->] (0,0) -- (1.5,0.5);
   \end{tikzpicture}};
  (40,10)*{
  \begin{tikzpicture} [scale=.4]
    \draw[very thick, black] (3.5,2.5) to [out=90,in=0] (3,3) to (1,3) to [out=180,in=90] (0.5,2.5);
    \draw[very thick, black] (3.5,0.5) to [out=270,in=0] (3,0) -- (1,0) to [out=180,in=270] (.5,0.5);
    \draw[very thick,black] (0.5,2.5) to [out=270,in=90] (1,1.5)  to  [out=270,in=90]  (.5,.5);
    \draw[very thick,black] (3.5,2.5) to [out=270,in=0] (3,2) to (2,2) to [out=180,in=90] (1.5,1.5) to [out=270,in=180] (2,1) to(3,1) to [out=0,in=90] (3.5,0.5);
  \end{tikzpicture}
  };
  (40,-10)*{
  \begin{tikzpicture} [scale=.4]
    \draw[very thick, black] (3.5,2.5) to [out=90,in=0] (3,3) to (1,3) to [out=180,in=90] (0.5,2.5);
    \draw[very thick, black] (3.5,0.5) to [out=270,in=0] (3,0) -- (1,0) to [out=180,in=270] (.5,0.5);
     \draw[very thick,black] (3.5,2.5) to [out=270,in=90] (3,1.5)  to  [out=270,in=90]  (3.5,.5);
    \draw[very thick,black] (0.5,2.5) to [out=270,in=180] (1,2) to (2,2) to [out=0,in=90] (2.5,1.5) to [out=270,in=0] (2,1) to(1,1) to [out=180,in=90] (0.5,0.5);
  \end{tikzpicture}
  };
  (5,-1)*{
  \uwave{\begin{tikzpicture} [scale=.4]
    \draw[very thick, black] (3.5,2.5) to [out=90,in=0] (3,3) to (1,3) to [out=180,in=90] (0.5,2.5);
    \draw[very thick, black] (3.5,0.5) to [out=270,in=0] (3,0) -- (1,0) to [out=180,in=270] (.5,0.5);
    \draw[very thick,black] (3.5,2.5) to [out=270,in=90] (3,1.5)  to  [out=270,in=90]  (3.5,.5);
    \draw[very thick,black] (.5,2.5) to [out=270,in=90] (1,1.5)  to  [out=270,in=90]  (.5,.5);
     \draw[very thick,black] (.5,2.5) to [out=270,in=90] (1,1.5)  to  [out=270,in=90]  (.5,.5);
     \draw[very thick] (2,1.5) circle (0.5);
  \end{tikzpicture}
  }
  };
  (32,10)*{q^{-1}};
  (32,-10)*{q^{-1}};
  (72,0)*{q^{-2}};
  (40,0)*{\bigoplus};
  \endxy 
  \end{equation}
  and can be interpreted as chain complex over the (additive) Bar-Natan
  category, alternatively in its $\Kh$- or $\mathrm{Lee}$-version. In the latter
  case, Bar-Natan--Morrison explain in \cite{MR2253455} how this chain complex splits upon
  proceeding to the Karoubi envelope. There, one uses colorings to encode images
  of the idempotents in $\mathrm{Lee}(\bigcirc) \cong \kb[X]/\langle X^2-1\rangle$. For each
  color, e.g. red, the direct summand
  \[
    \xy
    (80,0)*{
      \begin{tikzpicture} [scale=.4]
        \draw[very thick,dred] (3.5,2.5) to [out=90,in=0] (3,3) to (1,3) to [out=180,in=90] (0.5,2.5);
        \draw[very thick,dred] (3.5,0.5) to [out=270,in=0] (3,0) -- (1,0) to [out=180,in=270] (.5,0.5);
         \draw[very thick,dred] (0.5,2.5) to [out=270,in=180] (1,2) to (3,2) to [out=0,in=270] (3.5,2.5);
        \draw[very thick,dred]  (.5,0.5) to [out=90,in=180] (1,1) to (3,1) to [out=0,in=90] (3.5,0.5);
       \end{tikzpicture}
    };
    (88,0)*{
    \begin{tikzpicture} [scale=.4]
     \draw[opacity=0] (0,0) to (0,-5);
    \end{tikzpicture}
    };
    (20,5.5)*{\begin{tikzpicture} [scale=.5]
     \draw[->] (0,0) -- (1.5,0.5);
      \end{tikzpicture}};
     (20,-5.5)*{\begin{tikzpicture} [scale=.5]
     \draw[->] (0,0) -- (1.5,-0.5);
      \end{tikzpicture}};
     (60,5.5)*{\begin{tikzpicture} [scale=.5]
     \draw[->] (0,0) -- (1.5,-0.5);
      \end{tikzpicture}};
     (60,-5.5)*{\begin{tikzpicture} [scale=.5]
     \draw[->] (0,0) -- (1.5,0.5);
     \end{tikzpicture}};
    (40,10)*{
    \begin{tikzpicture} [scale=.4]
      \draw[very thick,dred] (3.5,2.5) to [out=90,in=0] (3,3) to (1,3) to [out=180,in=90] (0.5,2.5);
      \draw[very thick,dred] (3.5,0.5) to [out=270,in=0] (3,0) -- (1,0) to [out=180,in=270] (.5,0.5);
      \draw[very thick,dred] (0.5,2.5) to [out=270,in=90] (1,1.5)  to  [out=270,in=90]  (.5,.5);
      \draw[very thick,dred] (3.5,2.5) to [out=270,in=0] (3,2) to (2,2) to [out=180,in=90] (1.5,1.5) to [out=270,in=180] (2,1) to(3,1) to [out=0,in=90] (3.5,0.5);
    \end{tikzpicture}
    };
    (40,-10)*{
    \begin{tikzpicture} [scale=.4]
      \draw[very thick,dred] (3.5,2.5) to [out=90,in=0] (3,3) to (1,3) to [out=180,in=90] (0.5,2.5);
      \draw[very thick,dred] (3.5,0.5) to [out=270,in=0] (3,0) -- (1,0) to [out=180,in=270] (.5,0.5);
       \draw[very thick,dred] (3.5,2.5) to [out=270,in=90] (3,1.5)  to  [out=270,in=90]  (3.5,.5);
      \draw[very thick,dred] (0.5,2.5) to [out=270,in=180] (1,2) to (2,2) to [out=0,in=90] (2.5,1.5) to [out=270,in=0] (2,1) to(1,1) to [out=180,in=90] (0.5,0.5);
    \end{tikzpicture}
    };
    (5,-1)*{
    \uwave{\begin{tikzpicture} [scale=.4]
      \draw[very thick,dred] (3.5,2.5) to [out=90,in=0] (3,3) to (1,3) to [out=180,in=90] (0.5,2.5);
      \draw[very thick,dred] (3.5,0.5) to [out=270,in=0] (3,0) -- (1,0) to [out=180,in=270] (.5,0.5);
      \draw[very thick,dred] (3.5,2.5) to [out=270,in=90] (3,1.5)  to  [out=270,in=90]  (3.5,.5);
      \draw[very thick,dred] (.5,2.5) to [out=270,in=90] (1,1.5)  to  [out=270,in=90]  (.5,.5);
       \draw[very thick,dred] (.5,2.5) to [out=270,in=90] (1,1.5)  to  [out=270,in=90]  (.5,.5);
       \draw[very thick,dred] (2,1.5) circle (0.5);
    \end{tikzpicture}
    }
    };
    (32,10)*{q^{-1}};
    (32,-10)*{q^{-1}};
    (72,0)*{q^{-2}};
    (40,0)*{\bigoplus};
    \endxy 
  \]
  is contractible, as all components of the differential (induced by saddle
  cobordisms) are isomorphisms. The remaining components of the cube of
  resolutions in the Karoubi envelope are:
  \[
    \uwave{\begin{tikzpicture} [scale=.4,anchorbase]
      \draw[very thick,dred] (3.5,2.5) to [out=90,in=0] (3,3) to (1,3) to [out=180,in=90] (0.5,2.5);
      \draw[very thick,dred] (3.5,0.5) to [out=270,in=0] (3,0) -- (1,0) to [out=180,in=270] (.5,0.5);
      \draw[very thick,dred] (3.5,2.5) to [out=270,in=90] (3,1.5)  to  [out=270,in=90]  (3.5,.5);
      \draw[very thick,dred] (.5,2.5) to [out=270,in=90] (1,1.5)  to  [out=270,in=90]  (.5,.5);
       \draw[very thick,dred] (.5,2.5) to [out=270,in=90] (1,1.5)  to  [out=270,in=90]  (.5,.5);
       \draw[very thick,green] (2,1.5) circle (0.5);
    \end{tikzpicture}
    }
    \quad,\quad
    \uwave{\begin{tikzpicture} [scale=.4,anchorbase]
      \draw[very thick,green] (3.5,2.5) to [out=90,in=0] (3,3) to (1,3) to [out=180,in=90] (0.5,2.5);
      \draw[very thick,green] (3.5,0.5) to [out=270,in=0] (3,0) -- (1,0) to [out=180,in=270] (.5,0.5);
      \draw[very thick,green] (3.5,2.5) to [out=270,in=90] (3,1.5)  to  [out=270,in=90]  (3.5,.5);
      \draw[very thick,green] (.5,2.5) to [out=270,in=90] (1,1.5)  to  [out=270,in=90]  (.5,.5);
       \draw[very thick,green] (.5,2.5) to [out=270,in=90] (1,1.5)  to  [out=270,in=90]  (.5,.5);
       \draw[very thick,dred] (2,1.5) circle (0.5);
    \end{tikzpicture}
    }
    \quad,\quad
    q^{-2}t^2
    \begin{tikzpicture} [scale=.4,anchorbase]
      \draw[very thick,green] (3.5,2.5) to [out=90,in=0] (3,3) to (1,3) to [out=180,in=90] (0.5,2.5);
      \draw[very thick,dred] (3.5,0.5) to [out=270,in=0] (3,0) -- (1,0) to [out=180,in=270] (.5,0.5);
       \draw[very thick,green] (0.5,2.5) to [out=270,in=180] (1,2) to (3,2) to [out=0,in=270] (3.5,2.5);
      \draw[very thick,dred]  (.5,0.5) to [out=90,in=180] (1,1) to (3,1) to [out=0,in=90] (3.5,0.5);
     \end{tikzpicture}
     \quad,\quad
     q^{-2}t^2
    \begin{tikzpicture} [scale=-.4,anchorbase]
      \draw[very thick,green] (3.5,2.5) to [out=90,in=0] (3,3) to (1,3) to [out=180,in=90] (0.5,2.5);
      \draw[very thick,dred] (3.5,0.5) to [out=270,in=0] (3,0) -- (1,0) to [out=180,in=270] (.5,0.5);
       \draw[very thick,green] (0.5,2.5) to [out=270,in=180] (1,2) to (3,2) to [out=0,in=270] (3.5,2.5);
      \draw[very thick,dred]  (.5,0.5) to [out=90,in=180] (1,1) to (3,1) to [out=0,in=90] (3.5,0.5);
     \end{tikzpicture}
  \]
  These four components give rise to generators for the 4-dimensional Lee homology of the Hopf link. The components are in (non-canonical) bijection with orientations of the Hopf link. A
  preferred bijection can be determined after fixing a preferred root of the
  deformation polynomial $X^2-1$, say $1$.
\end{example}

\begin{example}\label{example:Hopfsplit} Now we work with the deformed foam
  category $\KhRS$ with $\Sigma=\{+1,-1\}$. In the Karoubi envelope, we again
  use colorings to encode images of the idempotents associated to $1$ (red, darker)
  and $-1$  (green, lighter). For complexes thereover we have the isomorphisms, \cite{MR3590355}:
\[
  \CC{\begin{tikzpicture}[anchorbase,scale=.2]
    \draw [very thick, green, ->] (1,-2)to(1,-1.7) to [out=90,in=270] (-1,1.7) to (-1,2);
    \draw [white, line width=.15cm] (-1,-2) to (-1,-1.7) to [out=90,in=270] (1,1.7) to (1,2);
    \draw [very thick,dred, ->] (-1,-2) to (-1,-1.7) to [out=90,in=270] (1,1.7) to (1,2);
  \end{tikzpicture}
  }
   \quad = \quad 
   \uwave{
  \begin{tikzpicture}[anchorbase,scale=.2]
    \draw [very thick,green, ->] (1,-2)to(1,-1.7) to [out=90,in=315] (0,-.5)(0,.5) to [out=135,in=270] (-1,1.7) to (-1,2);
    \draw[double] (0,.5) to (0,-.5);
    \draw [very thick,dred, ->] (-1,-2) to(-1,-1.7) to [out=90,in=225] (0,-.5)  (0,.5) to [out=45,in=270] (1,1.7) to (1,2);
  \end{tikzpicture}
  }
  \to
    q^{-1}\;
  \begin{tikzpicture}[anchorbase,scale=.2]
    \draw [very thick, green, ->] (1,-2) to (1,0) (-1,0) to (-1,2);
    \draw [very thick,dred, ->] (-1,-2) to (-1,0) (1,0) to (1,2);
  \end{tikzpicture}
  \quad \cong \quad
  \begin{tikzpicture}[anchorbase,scale=.2]
    \draw [very thick,green, ->] (1,-2)to(1,-1.7) to [out=90,in=315] (0,-.5)(0,.5) to [out=135,in=270] (-1,1.7) to (-1,2);
    \draw[double] (0,.5) to (0,-.5);
    \draw [very thick,dred, ->] (-1,-2) to(-1,-1.7) to [out=90,in=225] (0,-.5)  (0,.5) to [out=45,in=270] (1,1.7) to (1,2);
  \end{tikzpicture}
  \quad \cong \quad 
  q\;
  \begin{tikzpicture}[anchorbase,scale=.2]
    \draw [very thick, green, ->] (1,-2) to (1,0) (-1,0) to (-1,2);
    \draw [very thick,dred, ->] (-1,-2) to (-1,0) (1,0) to (1,2);
  \end{tikzpicture}
  \to
  \uwave{
  \begin{tikzpicture}[anchorbase,scale=.2]
    \draw [very thick,green, ->] (1,-2)to(1,-1.7) to [out=90,in=315] (0,-.5)(0,.5) to [out=135,in=270] (-1,1.7) to (-1,2);
    \draw[double] (0,.5) to (0,-.5);
    \draw [very thick,dred, ->] (-1,-2) to(-1,-1.7) to [out=90,in=225] (0,-.5)  (0,.5) to [out=45,in=270] (1,1.7) to (1,2);
  \end{tikzpicture}
  }
  \quad = \quad
  \CC{\begin{tikzpicture}[anchorbase,scale=.2]
    \draw [very thick,dred, ->] (-1,-2) to (-1,-1.7) to [out=90,in=270] (1,1.7) to (1,2);
    \draw [white, line width=.15cm] (1,-2)to(1,-1.7) to [out=90,in=270] (-1,1.7) to (-1,2);
    \draw [very thick,green, ->] (1,-2)to(1,-1.7) to [out=90,in=270] (-1,1.7) to (-1,2);
  \end{tikzpicture}
  }
\] 
Here we have used that webs with inconsistently colored edges are zero objects.
In particular, $1$-labelled strands with distinct colors braid trivially. For a
positive Hopf link with a coloring of its components by distinct colors, we thus
obtain:
 
  \begin{equation}
    \label{eqn:deformedHopf}  
  \xy
  (-18,0)*{
    \CC{\begin{tikzpicture} [xscale=-.4,yscale=.4,anchorbase]
   \draw[dred, very thick, directed=.6] (1.7,1.3) to [out=315,in=180] (2,1) to [out=0,in=180] (3,2) to [out=0,in=270] (3.5,2.5) to [out=90,in=0] (3,3) to (1,3) to [out=180,in=90] (0.5,2.5) to [out=270,in=180] (1,2) to [out=0,in=135] (1.3,1.7);
   \draw[green, very thick,directed=.4] (2.7,1.3) to [out=315,in=180] (3,1) to [out=0,in=90] (3.5,0.5) to [out=270,in=0] (3,0) -- (1,0) to [out=180,in=270] (.5,0.5) to [out=90,in=180] (1,1) to [out=0,in=180] (2,2) to [out=0,in=135] (2.3,1.7);
  \end{tikzpicture}} =};
  (80,0)*{
  \begin{tikzpicture} [scale=.4]
   \draw[very thick, dred] (3.5,2.5) to [out=90,in=0] (3,3) to (1,3) to [out=180,in=90] (0.5,2.5);
   \draw[very thick, green] (3.5,0.5) to [out=270,in=0] (3,0) -- (1,0) to [out=180,in=270] (.5,0.5);
    \draw[very thick, dred] (0.5,2.5) to [out=270,in=180] (1,2) (3,2) to [out=0,in=270] (3.5,2.5) (1,1) to (3,1);
   \draw[very thick,green]  (.5,0.5) to [out=90,in=180] (1,1) (3,1) to [out=0,in=90] (3.5,0.5) (1,2) to (3,2);
  \end{tikzpicture}
  };
  (88,0)*{
  \begin{tikzpicture} [scale=.4]
   \draw[opacity=0] (0,0) to (0,-5);
  \end{tikzpicture}
  };
  (20,5.5)*{\begin{tikzpicture} [scale=.5]
   \draw[->] (0,0) -- (1.5,0.5);
    \end{tikzpicture}};
   (20,-5.5)*{\begin{tikzpicture} [scale=.5]
   \draw[->] (0,0) -- (1.5,-0.5);
    \end{tikzpicture}};
   (60,5.5)*{\begin{tikzpicture} [scale=.5]
   \draw[->] (0,0) -- (1.5,-0.5);
    \end{tikzpicture}};
   (60,-5.5)*{\begin{tikzpicture} [scale=.5]
   \draw[->] (0,0) -- (1.5,0.5);
   \end{tikzpicture}};
  (40,10)*{
  \begin{tikzpicture} [scale=.4]
    \draw[very thick, dred] (3,2) \ru (3.5,2.5) to [out=90,in=0] (3,3) to (1,3) to [out=180,in=90] (0.5,2.5) (1.5,1.5)\dr (2,1) to (3,1);
    \draw[very thick, green] (3,1) \rd (3.5,0.5) to [out=270,in=0] (3,0) -- (1,0) to [out=180,in=270] (.5,0.5)  (1.5,1.5)\ur (2,2) to (3,2);
   \draw[double] (1,1.5) to (1.5,1.5);
   \draw[very thick,dred] (.5,2.5) to [out=270,in=90] (1,1.5);
   \draw[very thick,green] (1,1.5)  to  [out=270,in=90]  (.5,.5);
  \end{tikzpicture}
  };
  (40,-10)*{
  \begin{tikzpicture} [scale=.4]
    \draw[very thick, dred] (3.5,2.5) to [out=90,in=0] (3,3) to (1,3) to [out=180,in=90] (0.5,2.5) \dr (1,2) (2.5,1.5) \dl (2,1) to (1,1);
    \draw[very thick, green] (3.5,0.5) to [out=270,in=0] (3,0) -- (1,0) to [out=180,in=270] (.5,0.5) \ur (1,1) (2.5,1.5) \ul (2,2) to (1,2);
  \draw[double] (2.5,1.5) to (3,1.5);
  \draw[very thick, dred] (3.5,2.5) to [out=270,in=90] (3,1.5);
  \draw[very thick,green] (3,1.5)  to  [out=270,in=90] (3.5,.5);
  \end{tikzpicture}
  };
  (5,-1)*{
  \uwave{\begin{tikzpicture} [scale=.4]
    \draw[very thick, dred] (3.5,2.5) to [out=90,in=0] (3,3) to (1,3) to [out=180,in=90] (0.5,2.5);
    \draw[very thick, green] (3.5,0.5) to [out=270,in=0] (3,0) -- (1,0) to [out=180,in=270] (.5,0.5);
   \draw[very thick, dred] (3.5,2.5) to [out=270,in=90] (3,1.5) (.5,2.5) to [out=270,in=90] (1,1.5);
   \draw[very thick,green] (1,1.5)  to  [out=270,in=90]  (.5,.5) (3,1.5)  to  [out=270,in=90] (3.5,.5);
    \draw[very thick,dred] (1.5,1.5) \dr (2,1) \ru (2.5,1.5);
    \draw[very thick,green] (1.5,1.5) \ur (2,2) \rd (2.5,1.5);
    \draw[double] (1,1.5) to (1.5,1.5);
    \draw[double] (2.5,1.5) to (3,1.5);
  \end{tikzpicture}
  }
  };
  (32,10)*{q^{-1}};
  (32,-10)*{q^{-1}};
  (72,0)*{q^{-2}};
  (40,0)*{\bigoplus};
  \endxy 
  =
  \uwave{\begin{tikzpicture} [scale=.4,anchorbase]
    \draw[very thick, dred] (3.5,2.5) to [out=90,in=0] (3,3) to (1,3) to [out=180,in=90] (0.5,2.5);
    \draw[very thick, green] (3.5,0.5) to [out=270,in=0] (3,0) -- (1,0) to [out=180,in=270] (.5,0.5);
   \draw[very thick, dred] (3.5,2.5) to [out=270,in=90] (3,1.5) (.5,2.5) to [out=270,in=90] (1,1.5);
   \draw[very thick,green] (1,1.5)  to  [out=270,in=90]  (.5,.5) (3,1.5)  to  [out=270,in=90] (3.5,.5);
    \draw[very thick,dred] (1.5,1.5) \dr (2,1) \ru (2.5,1.5);
    \draw[very thick,green] (1.5,1.5) \ur (2,2) \rd (2.5,1.5);
    \draw[double] (1,1.5) to (1.5,1.5);
    \draw[double] (2.5,1.5) to (3,1.5);
  \end{tikzpicture}
  }
  \cong 
  \uwave{\begin{tikzpicture} [scale=.4,anchorbase]
    \draw[very thick, dred] (3.5,2.5) to [out=90,in=0] (3,3) to (1,3) to [out=180,in=90] (0.5,2.5);
    \draw[very thick, green] (3.5,0.5) to [out=270,in=0] (3,0) -- (1,0) to [out=180,in=270] (.5,0.5);
     \draw[very thick, dred] (0.5,2.5) to [out=270,in=180] (1,2) (3,2) to [out=0,in=270] (3.5,2.5) (1,2) to (3,2);
    \draw[very thick,green]  (.5,0.5) to [out=90,in=180] (1,1) (3,1) to [out=0,in=90] (3.5,0.5) (1,1) to (3,1);
   \end{tikzpicture}
  }
  \end{equation}
For equal colors, on the other hand, we obtain
\[
  \CC{\begin{tikzpicture}[anchorbase,scale=.2]
    \draw [very thick, dred, ->] (1,-2)to(1,-1.7) to [out=90,in=270] (-1,1.7) to (-1,2);
    \draw [white, line width=.15cm] (-1,-2) to (-1,-1.7) to [out=90,in=270] (1,1.7) to (1,2);
    \draw [very thick,dred, ->] (-1,-2) to (-1,-1.7) to [out=90,in=270] (1,1.7) to (1,2);
  \end{tikzpicture}
  }
   \quad = \quad 
   \uwave{
  \begin{tikzpicture}[anchorbase,scale=.2]
    \draw [very thick,dred, ->] (1,-2)to(1,-1.7) to [out=90,in=315] (0,-.5)(0,.5) to [out=135,in=270] (-1,1.7) to (-1,2);
    \draw[double] (0,.5) to (0,-.5);
    \draw [very thick,dred, ->] (-1,-2) to(-1,-1.7) to [out=90,in=225] (0,-.5)  (0,.5) to [out=45,in=270] (1,1.7) to (1,2);
  \end{tikzpicture}
  }
  \to
    q^{-1}\;
  \begin{tikzpicture}[anchorbase,scale=.2]
    \draw [very thick, dred, ->] (1,-2) to (1,0) (-1,0) to (-1,2);
    \draw [very thick,dred, ->] (-1,-2) to (-1,0) (1,0) to (1,2);
  \end{tikzpicture}
  \quad \cong \quad
  q^{-1} t\;
  \begin{tikzpicture}[anchorbase,scale=.2]
    \draw [very thick, dred, ->] (1,-2) to (1,0) (-1,0) to (-1,2);
    \draw [very thick,dred, ->] (-1,-2) to (-1,0) (1,0) to (1,2);
  \end{tikzpicture}
  \qquad,\qquad 
  \CC{\begin{tikzpicture}[anchorbase,scale=.2]
    \draw [very thick,dred, ->] (-1,-2) to (-1,-1.7) to [out=90,in=270] (1,1.7) to (1,2);
    \draw [white, line width=.15cm] (1,-2)to(1,-1.7) to [out=90,in=270] (-1,1.7) to (-1,2);
    \draw [very thick,dred, ->] (1,-2)to(1,-1.7) to [out=90,in=270] (-1,1.7) to (-1,2);
  \end{tikzpicture}
  }
  \quad = \quad
  q\;
  \begin{tikzpicture}[anchorbase,scale=.2]
    \draw [very thick, dred, ->] (1,-2) to (1,0) (-1,0) to (-1,2);
    \draw [very thick,dred, ->] (-1,-2) to (-1,0) (1,0) to (1,2);
  \end{tikzpicture}
  \to
  \uwave{
  \begin{tikzpicture}[anchorbase,scale=.2]
    \draw [very thick,dred, ->] (1,-2)to(1,-1.7) to [out=90,in=315] (0,-.5)(0,.5) to [out=135,in=270] (-1,1.7) to (-1,2);
    \draw[double] (0,.5) to (0,-.5);
    \draw [very thick,dred, ->] (-1,-2) to(-1,-1.7) to [out=90,in=225] (0,-.5)  (0,.5) to [out=45,in=270] (1,1.7) to (1,2);
  \end{tikzpicture}
  }
  \quad \cong \quad
  q t^{-1}\;
  \begin{tikzpicture}[anchorbase,scale=.2]
    \draw [very thick, dred, ->] (1,-2) to (1,0) (-1,0) to (-1,2);
    \draw [very thick,dred, ->] (-1,-2) to (-1,0) (1,0) to (1,2);
  \end{tikzpicture}
\] 
since the label $2$ on the double edge exceeds the multiplicity of the
corresponding element (here $1$) in $\Sigma$. For the consistently-colored Hopf link, we thus obtain:

\begin{equation}
 \label{eqn:deformedHopftwo}  
\xy
(-18,0)*{
  \CC{\begin{tikzpicture} [xscale=-.4,yscale=.4,anchorbase]
 \draw[dred, very thick, directed=.6] (1.7,1.3) to [out=315,in=180] (2,1) to [out=0,in=180] (3,2) to [out=0,in=270] (3.5,2.5) to [out=90,in=0] (3,3) to (1,3) to [out=180,in=90] (0.5,2.5) to [out=270,in=180] (1,2) to [out=0,in=135] (1.3,1.7);
 \draw[dred, very thick,directed=.4] (2.7,1.3) to [out=315,in=180] (3,1) to [out=0,in=90] (3.5,0.5) to [out=270,in=0] (3,0) -- (1,0) to [out=180,in=270] (.5,0.5) to [out=90,in=180] (1,1) to [out=0,in=180] (2,2) to [out=0,in=135] (2.3,1.7);
\end{tikzpicture}} =};
(80,0)*{
\begin{tikzpicture} [scale=.4]
 \draw[very thick, dred] (3.5,2.5) to [out=90,in=0] (3,3) to (1,3) to [out=180,in=90] (0.5,2.5);
 \draw[very thick, dred] (3.5,0.5) to [out=270,in=0] (3,0) -- (1,0) to [out=180,in=270] (.5,0.5);
  \draw[very thick, dred] (0.5,2.5) to [out=270,in=180] (1,2) (3,2) to [out=0,in=270] (3.5,2.5) (1,1) to (3,1);
 \draw[very thick,dred]  (.5,0.5) to [out=90,in=180] (1,1) (3,1) to [out=0,in=90] (3.5,0.5) (1,2) to (3,2);
\end{tikzpicture}
};
(88,0)*{
\begin{tikzpicture} [scale=.4]
 \draw[opacity=0] (0,0) to (0,-5);
\end{tikzpicture}
};
(20,5.5)*{\begin{tikzpicture} [scale=.5]
 \draw[->] (0,0) -- (1.5,0.5);
  \end{tikzpicture}};
 (20,-5.5)*{\begin{tikzpicture} [scale=.5]
 \draw[->] (0,0) -- (1.5,-0.5);
  \end{tikzpicture}};
 (60,5.5)*{\begin{tikzpicture} [scale=.5]
 \draw[->] (0,0) -- (1.5,-0.5);
  \end{tikzpicture}};
 (60,-5.5)*{\begin{tikzpicture} [scale=.5]
 \draw[->] (0,0) -- (1.5,0.5);
 \end{tikzpicture}};
(40,10)*{
\begin{tikzpicture} [scale=.4]
  \draw[very thick, dred] (3,2) \ru (3.5,2.5) to [out=90,in=0] (3,3) to (1,3) to [out=180,in=90] (0.5,2.5) (1.5,1.5)\dr (2,1) to (3,1);
  \draw[very thick, dred] (3,1) \rd (3.5,0.5) to [out=270,in=0] (3,0) -- (1,0) to [out=180,in=270] (.5,0.5)  (1.5,1.5)\ur (2,2) to (3,2);
 \draw[double] (1,1.5) to (1.5,1.5);
 \draw[very thick,dred] (.5,2.5) to [out=270,in=90] (1,1.5);
 \draw[very thick,dred] (1,1.5)  to  [out=270,in=90]  (.5,.5);
\end{tikzpicture}
};
(40,-10)*{
\begin{tikzpicture} [scale=.4]
  \draw[very thick, dred] (3.5,2.5) to [out=90,in=0] (3,3) to (1,3) to [out=180,in=90] (0.5,2.5) \dr (1,2) (2.5,1.5) \dl (2,1) to (1,1);
  \draw[very thick, dred] (3.5,0.5) to [out=270,in=0] (3,0) -- (1,0) to [out=180,in=270] (.5,0.5) \ur (1,1) (2.5,1.5) \ul (2,2) to (1,2);
\draw[double] (2.5,1.5) to (3,1.5);
\draw[very thick, dred] (3.5,2.5) to [out=270,in=90] (3,1.5);
\draw[very thick,dred] (3,1.5)  to  [out=270,in=90] (3.5,.5);
\end{tikzpicture}
};
(5,-1)*{
\uwave{\begin{tikzpicture} [scale=.4]
  \draw[very thick, dred] (3.5,2.5) to [out=90,in=0] (3,3) to (1,3) to [out=180,in=90] (0.5,2.5);
  \draw[very thick, dred] (3.5,0.5) to [out=270,in=0] (3,0) -- (1,0) to [out=180,in=270] (.5,0.5);
 \draw[very thick, dred] (3.5,2.5) to [out=270,in=90] (3,1.5) (.5,2.5) to [out=270,in=90] (1,1.5);
 \draw[very thick,dred] (1,1.5)  to  [out=270,in=90]  (.5,.5) (3,1.5)  to  [out=270,in=90] (3.5,.5);
  \draw[very thick,dred] (1.5,1.5) \dr (2,1) \ru (2.5,1.5);
  \draw[very thick,dred] (1.5,1.5) \ur (2,2) \rd (2.5,1.5);
  \draw[double] (1,1.5) to (1.5,1.5);
  \draw[double] (2.5,1.5) to (3,1.5);
\end{tikzpicture}
}
};
(32,10)*{q^{-1}};
(32,-10)*{q^{-1}};
(72,0)*{q^{-2}};
(40,0)*{\bigoplus};
\endxy 
=
q^{-2}t^2\uwave{\begin{tikzpicture} [scale=.4,anchorbase]
  \draw[very thick, dred] (3.5,2.5) to [out=90,in=0] (3,3) to (1,3) to [out=180,in=90] (0.5,2.5);
  \draw[very thick, dred] (3.5,0.5) to [out=270,in=0] (3,0) -- (1,0) to [out=180,in=270] (.5,0.5);
   \draw[very thick, dred] (0.5,2.5) to [out=270,in=180] (1,2) (3,2) to [out=0,in=270] (3.5,2.5) (1,2) to (3,2);
  \draw[very thick,dred]  (.5,0.5) to [out=90,in=180] (1,1) (3,1) to [out=0,in=90] (3.5,0.5) (1,1) to (3,1);
 \end{tikzpicture}
}
\end{equation}
These colorings, together with the symmetric configurations with red and green
swapped, generate the four non-trivial components of the Hopf link complex and
give rise to the four components in the $\glnn{2}$-Lee homology of the Hopf
link. It is clear that these summands are naturally parametrized by colorings of
the link components by elements of $\Sigma$. 

Conversely, one can also consider all admissible colorings of the webs that
appear in the cube of resolutions of the Hopf link complex. The complex then
splits into a direct sum of the complexes already observed above, together with:

\begin{equation}
  \xy
  (80,0)*{
  \begin{tikzpicture} [scale=.4]
   \draw[very thick, dred] (3.5,2.5) to [out=90,in=0] (3,3) to (1,3) to [out=180,in=90] (0.5,2.5);
   \draw[very thick, green] (3.5,0.5) to [out=270,in=0] (3,0) -- (1,0) to [out=180,in=270] (.5,0.5);
    \draw[very thick, dred] (0.5,2.5) to [out=270,in=180] (1,2) (3,2) to [out=0,in=270] (3.5,2.5) (1,2) to (3,2);
   \draw[very thick,green]  (.5,0.5) to [out=90,in=180] (1,1) (3,1) to [out=0,in=90] (3.5,0.5) (1,1) to (3,1);
  \end{tikzpicture}
  };
  (88,0)*{
  \begin{tikzpicture} [scale=.4]
   \draw[opacity=0] (0,0) to (0,-5);
  \end{tikzpicture}
  };
  (20,5.5)*{\begin{tikzpicture} [scale=.5]
   \draw[->] (0,0) -- (1.5,0.5);
    \end{tikzpicture}};
   (20,-5.5)*{\begin{tikzpicture} [scale=.5]
   \draw[->] (0,0) -- (1.5,-0.5);
    \end{tikzpicture}};
   (60,5.5)*{\begin{tikzpicture} [scale=.5]
   \draw[->] (0,0) -- (1.5,-0.5);
    \end{tikzpicture}};
   (60,-5.5)*{\begin{tikzpicture} [scale=.5]
   \draw[->] (0,0) -- (1.5,0.5);
   \end{tikzpicture}};
  (40,10)*{
  \begin{tikzpicture} [scale=.4]
    \draw[very thick, dred] (3,2) \ru (3.5,2.5) to [out=90,in=0] (3,3) to (1,3) to [out=180,in=90] (0.5,2.5) (1.5,1.5)\ur (2,2) to (3,2);
    \draw[very thick, green] (3,1) \rd (3.5,0.5) to [out=270,in=0] (3,0) -- (1,0) to [out=180,in=270] (.5,0.5) (1.5,1.5)\dr (2,1) to (3,1);
   \draw[double] (1,1.5) to (1.5,1.5);
   \draw[very thick,dred] (.5,2.5) to [out=270,in=90] (1,1.5);
   \draw[very thick,green] (1,1.5)  to  [out=270,in=90]  (.5,.5);
  \end{tikzpicture}
  };
  (40,-10)*{
  \begin{tikzpicture} [scale=.4]
    \draw[very thick, dred] (3.5,2.5) to [out=90,in=0] (3,3) to (1,3) to [out=180,in=90] (0.5,2.5) \dr (1,2) (2.5,1.5) \ul (2,2) to (1,2);
    \draw[very thick, green] (3.5,0.5) to [out=270,in=0] (3,0) -- (1,0) to [out=180,in=270] (.5,0.5) \ur (1,1) (2.5,1.5) \dl (2,1) to (1,1);
  \draw[double] (2.5,1.5) to (3,1.5);
  \draw[very thick, dred] (3.5,2.5) to [out=270,in=90] (3,1.5);
  \draw[very thick,green] (3,1.5)  to  [out=270,in=90] (3.5,.5);
  \end{tikzpicture}
  };
  (5,-1)*{
  \uwave{\begin{tikzpicture} [scale=.4]
    \draw[very thick, dred] (3.5,2.5) to [out=90,in=0] (3,3) to (1,3) to [out=180,in=90] (0.5,2.5);
    \draw[very thick, green] (3.5,0.5) to [out=270,in=0] (3,0) -- (1,0) to [out=180,in=270] (.5,0.5);
   \draw[very thick, dred] (3.5,2.5) to [out=270,in=90] (3,1.5) (.5,2.5) to [out=270,in=90] (1,1.5);
   \draw[very thick,green] (1,1.5)  to  [out=270,in=90]  (.5,.5) (3,1.5)  to  [out=270,in=90] (3.5,.5);
    \draw[very thick,green] (1.5,1.5) \dr (2,1) \ru (2.5,1.5);
    \draw[very thick,dred] (1.5,1.5) \ur (2,2) \rd (2.5,1.5);
    \draw[double] (1,1.5) to (1.5,1.5);
    \draw[double] (2.5,1.5) to (3,1.5);
  \end{tikzpicture}
  }
  };
  (32,10)*{q^{-1}};
  (32,-10)*{q^{-1}};
  (72,0)*{q^{-2}};
  (40,0)*{\bigoplus};
  \endxy 
  \end{equation}  
and a symmetrically colored version, both of which are contractible, as all components of the differential are isomorphisms.
  \end{example}

  \begin{example}
  As in Example~\ref{example:Hopfsplit}, one can compute deformed homology
  $\KhRS$ of the Hopf link for $\Sigma=\{\lambda_1^M, \lambda_2^{N-M}\}$. The
  computation \eqref{eqn:deformedHopf} generalizes directly and identifies a
  summand isomorphic to a tensor product of $\glnn{M}$ and $\glnn{N-M}$
  homologies of unknots. The consistently colored summands directly compute $\glnn{M}$ and $\glnn{N-M}$ homologies of the Hopf link. Note that the last equality in \eqref{eqn:deformedHopftwo} only holds in the $M=N=1$ case.
  \end{example}

\section{Functoriality considerations}
The original construction of Khovanov homology~\cite{MR1740682} associates chain
maps to link cobordisms in $\R^3\times I$ (represented as movies of link
diagrams) and Khovanov conjectured that these are independent under isotopy up to
homotopy and possibly factors of $-1$. This conjecture was verified by
Jacobsson~\cite{MR2113903}, Bar-Natan~\cite{MR2174270} and
Khovanov~\cite{MR2171235}, however the sign-issue proved persistent in the
original construction (but see Remark~\ref{rem:cobKh}). 

\subsection{Symmetry breaking for the crossing}
Modified constructions of Khovanov homology that fix the sign ambiguity in
the functoriality under link cobordisms were subsequently developed by
Caprau~\cite{MR2443094}, Clark--Morrison--Walker~\cite{MR2496052},
Blanchet~\cite{MR2647055}, Vogel~\cite{MR4096813}. What these approaches have in
common is that they use the orientation of the link diagram in a more
significant way than Khovanov homology to break the symmetry of the basic chain
complex associated to a single crossing. To illustrate this difference, recall
from \cite{MR2174270} that the Khovanov complexes for positive and negative
crossings

\begin{equation*}
  \CC{
  \begin{tikzpicture}[anchorbase,scale=.45]
    \draw [very thick, ->] (1,0) \pu (0,1.5);
    \draw [white, line width=.15cm] (0,0) \pu (1,1.5);
    \draw [very thick, ->] (0,0) \pu (1,1.5);
  \end{tikzpicture}
  }
  :=\left( 0 \to 
  \uwave{q^{1} \CC{
      \begin{tikzpicture}[anchorbase,scale=.45]
        \draw [very thick] (0,0) -- (0,1.5);
        \draw [very thick] (1,0) -- (1,1.5);
      \end{tikzpicture}
    }} \to q^{2}
    \CC{
      \begin{tikzpicture}[anchorbase,scale=.45]
        \draw [very thick] (0,0) \ur (.5,0.5) \rd (1,0);
        \draw [very thick] (0,1.5) \dr (.5,1) \ru (1,1.5);
      \end{tikzpicture} 
    } \to 0
    \right)
    ,\;\;
    \CC{
  \begin{tikzpicture}[anchorbase,xscale=-.5,yscale=.5]
    \draw [very thick, ->] (1,0) \pu (0,1.5);
    \draw [white, line width=.15cm] (0,0) \pu (1,1.5);
    \draw [very thick, ->] (0,0) \pu (1,1.5);
  \end{tikzpicture}
  }
  :=\left( 0 \to q^{-2}
  \CC{
  \begin{tikzpicture}[anchorbase,scale=.45]
    \draw [very thick] (0,0) \ur (.5,0.5) \rd (1,0);
    \draw [very thick] (0,1.5) \dr (.5,1) \ru (1,1.5);
  \end{tikzpicture} 
  } \to 
  \uwave{q^{-1}\CC{
\begin{tikzpicture}[anchorbase,scale=.45]
  \draw [very thick] (0,0) -- (0,1.5);
  \draw [very thick] (1,0) -- (1,1.5);
\end{tikzpicture}
  }} \to 0
  \right)
\end{equation*}   
 only differ by a 90 degree rotation and a grading shift, while Blanchet's
 construction \cite{MR2647055} uses    
 \begin{equation*}
  \CC{
  \begin{tikzpicture}[anchorbase,scale=.45]
    \draw [very thick, ->] (1,0) \pu (0,1.5);
    \draw [white, line width=.15cm] (0,0) \pu (1,1.5);
    \draw [very thick, ->] (0,0) \pu (1,1.5);
  \end{tikzpicture}
  }
  :=\left( 0 \to 
  \uwave{q^{-1}\CC{
\begin{tikzpicture}[anchorbase,scale=.45]
\draw [very thick, ->] (0,0) -- (0,1.5);
\draw [very thick, ->] (1,0) -- (1,1.5);
\end{tikzpicture}
}} \to q^{-2}
\CC{
\begin{tikzpicture}[anchorbase,xscale=-.5,yscale=.5]
\draw [very thick] (0,0) \pu (0,0.1) \pu (.5,.6);
\draw [very thick] (1,0) \pu (1,0.1) \pu (.5,.6);
\draw [very thick, directed=.7] (.5,.6) \pu (.5,.9);
\node at (.8,.75) {\tiny $2$};
\draw [very thick,->] (.5,.9) \pu (0,1.3) \pu (0,1.5);
\draw [very thick,->] (.5,.9) \pu (1,1.3) \pu (1,1.5);
\end{tikzpicture} 
}
\to 0
\right) 
  ,   \;\;
    \CC{
  \begin{tikzpicture}[anchorbase,xscale=-.5,yscale=.5]
    \draw [very thick, ->] (1,0) \pu (0,1.5);
    \draw [white, line width=.15cm] (0,0) \pu (1,1.5);
    \draw [very thick, ->] (0,0) \pu (1,1.5);
  \end{tikzpicture}
  }
  :=\left( 0 \to q^{2}
    \CC{
    \begin{tikzpicture}[anchorbase,scale=.45]
      \draw [very thick] (0,0) \pu (0,0.1) \pu (.5,.6);
      \draw [very thick] (1,0) \pu (1,0.1) \pu (.5,.6);
      \draw [very thick, directed=.7] (.5,.6) \pu (.5,.9);
      \node at (.8,.75) {\tiny $2$};
      \draw [very thick,->] (.5,.9) \pu (0,1.3) \pu (0,1.5);
      \draw [very thick,->] (.5,.9) \pu (1,1.3) \pu (1,1.5);
    \end{tikzpicture} 
    } \to 
    \uwave{ q^{1}\CC{
  \begin{tikzpicture}[anchorbase,scale=.45]
    \draw [very thick, ->] (0,0) -- (0,1.5);
    \draw [very thick, ->] (1,0) -- (1,1.5);
  \end{tikzpicture}
    }} \to 0
    \right)
\end{equation*}  
where the symmetry is broken. Instead of the cup-cap pattern in Khovanov's
complex, an object associated to a 2-labelled edge is used. As argued in
Section~\ref{sec:skeinrel}, this can be though of as encoding the determinant
representation $\wedge^{2}(V)$, which is non-trivial for $\gltw$. In this sense
Blanchet's construction makes Khovanov homology functorial by leveraging
precisely the additionl data encoded in the diagrammatic calculus for the
representation category of $\Uqgltw$.

\begin{remark}[On functorial cobordism maps in the original Khovanov homology]
  \label{rem:cobKh}
The modified constructions of Khovanov homology mentioned above produce link
invariants that are (not necessarily canonically) isomorphic to the original
Khovanov homology. Having chosen such isomorphisms for all links, one can
transport the definition of functorial cobordism maps into the original theory.
The challenge is then to give an intrinsic description of such sign-corrected
maps, see Ehrig--Stroppel--Tubbenhauer~\cite{2016arXiv160108010E} and
Beliakova--Hogancamp--Putyra--Wehrli~\cite{2019arXiv190312194B}. An alternative
approach to fixing signs in Khovanov homology uses that the same deficiency
appears in Lee homology, where it is more straightforward to remedy, see also
Grigsby--Licata--Wehrli~\cite[Section 7.2]{MR3731256} and
Sano~\cite{2020arXiv200802131S}.
\end{remark}

Note that Blanchet's complexes, shown above, agree with the pattern
\eqref{eq:crossingcx} for the mirrored tangle diagram with the uniquely determined global shift that makes the homotopy equivalence realizing the Reidemeister 1 move in a 1-colored strand grading-preserving. And indeed, Blanchet's
version of Khovanov homology naturally fits into the family of
Khovanov--Rozansky homologies constructed via Robert--Wagner foams for $\glN$
\cite{1702.04140, MR3877770}. On the other hand, there only exists one example
of a link homology for $N>2$ that has a distinct $\slN$ flavour, namely
Khovanov's $\slnn{3}$ homology from \cite{MR2100691}.

\subsection{Base point action}
Any functorial link homology theory $H$ (and also some not properly functorial
constructions) assigns an algebra $H(\bigcirc)$ to the standard (oriented,
uncolored) unknot diagram, with multiplication induced by the
standardly-embedded pair of pants morphisms. Similarly, every base point on a
link diagram $L$ (uncolored) exhibits $H(L)$ as an $H(\bigcirc)$-module. If
$H(\bigcirc)$ is a quotient of a polynomial ring generated by $X$ in degree $2$,
such as in Khovanov and Khovanov--Rozansky homology, then it is enough to
understand the action of $X$ a.k.a. \emph{the dot}. 

In Khovanov homology as constructed via Blanchet foams and more generally in
Khovanov--Rozansky homology, these module structures only depend on the chosen
link component, not on the exact location of the base point. This is because the
action of $X$ one one side of a crossing is homotopic to action on the other
side:

\begin{equation}
  \label{eq:dotcrossing}
  \CC{
  \begin{tikzpicture}[anchorbase,scale=.45]
    \draw [very thick, ->] (1,0) \pu (0,1.5);
    \draw [white, line width=.15cm] (0,0) \pu (1,1.5);
    \draw [very thick, ->] (0,0) \pu (1,1.5);
    \node[opacity=1] at (.85,1) {$\bullet$};
  \end{tikzpicture}
  \!
  }  
  \simeq 
  \CC{
    \!
  \begin{tikzpicture}[anchorbase,scale=.45]
    \draw [very thick, ->] (1,0) \pu (0,1.5);
    \draw [white, line width=.15cm] (0,0) \pu (1,1.5);
    \draw [very thick, ->] (0,0) \pu (1,1.5);
    \node[opacity=1] at (.15,.35) {$\bullet$};
  \end{tikzpicture}
  }
  ,\qquad 
  \CC{
    \begin{tikzpicture}[anchorbase,scale=.45]
      \draw [very thick, ->] (1,0) \pu (0,1.5);
      \draw [white, line width=.15cm] (0,0) \pu (1,1.5);
      \draw [very thick, ->] (0,0) \pu (1,1.5);
      \node[opacity=1] at (.15,1) {$\bullet$};
    \end{tikzpicture}
    \!
    }  
    \simeq 
    \CC{
      \!
    \begin{tikzpicture}[anchorbase,scale=.45]
      \draw [very thick, ->] (1,0) \pu (0,1.5);
      \draw [white, line width=.15cm] (0,0) \pu (1,1.5);
      \draw [very thick, ->] (0,0) \pu (1,1.5);
      \node[opacity=1] at (.85,.35) {$\bullet$};
    \end{tikzpicture}
    }  
  \end{equation} 
  and similarly for negative crossings. For example, in Blanchet's complex for the positive crossing
\begin{equation*}
  \CC{
    \begin{tikzpicture}[anchorbase,scale=.45]
      \draw [very thick, ->] (1,0) \pu (0,1.5);
      \draw [white, line width=.15cm] (0,0) \pu (1,1.5);
      \draw [very thick, ->] (0,0) \pu (1,1.5);
    \end{tikzpicture}
    }
    :=\left( 0 \to 
    \uwave{q^{-1}\CC{
  \begin{tikzpicture}[anchorbase,scale=.45]
  \draw [very thick, ->] (0,0) -- (0,1.5);
  \draw [very thick, ->] (1,0) -- (1,1.5);
  \end{tikzpicture}
  }} \to q^{-2}
  \CC{
  \begin{tikzpicture}[anchorbase,xscale=-.5,yscale=.5]
  \draw [very thick] (0,0) \pu (0,0.1) \pu (.5,.6);
  \draw [very thick] (1,0) \pu (1,0.1) \pu (.5,.6);
  \draw [very thick, directed=.7] (.5,.6) \pu (.5,.9);
  \node at (.8,.75) {\tiny $2$};
  \draw [very thick,->] (.5,.9) \pu (0,1.3) \pu (0,1.5);
  \draw [very thick,->] (.5,.9) \pu (1,1.3) \pu (1,1.5);
  \end{tikzpicture} 
  }
  \to 0
  \right) 
\end{equation*}
the differential is given by a zip foam. Now consider the homotopy given by the
unzip foam pointing in the opposite direction. This is a homotopy for the
difference of $X$ actions on both sides of the crossing (either case of
\eqref{eq:dotcrossing}), as can be seen using the Blanchet foam relations. For
example, in degree zero:
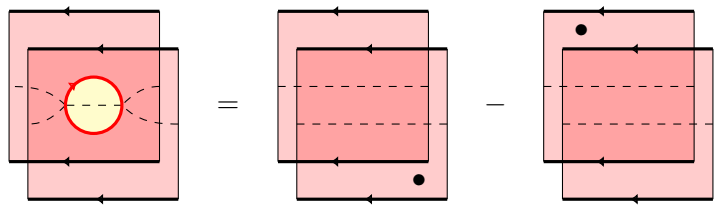
\[\begin{tikzpicture} [anchorbase,scale=.5,fill opacity=0.2]
\path[fill=red]  (2,5) to [out=90,in=270] (2,1) to [out=180,in=0] (0.25,1) to  (0.25,1.75) to [out=0, in=270] (1,2.5) to [out=90, in=0] (0.25,3.25) to [out=180, in=90] (-0.5,2.5) to [out=270,in=180] (0.25,1.75) to (0.25,1)  to (-2,1) to (-2,5) to (2,5);
\path[fill=red]  (2.5,4) to [out=90,in=270] (2.5,0) to [out=180,in=0] (0.25,0) to  (0.25,1.75) to [out=0, in=270] (1,2.5) to [out=90, in=0] (0.25,3.25) to [out=180, in=90] (-0.5,2.5) to [out=270,in=180] (0.25,1.75) to (0.25,0)  to (-1.5,0) to (-1.5,4) to (2.5,4);
\path[fill=yellow] (0.25,1.75) to [out=0, in=270] (1,2.5) to [out=90, in=0] (0.25,3.25) to [out=180, in=90] (-0.5,2.5) to [out=270,in=180] (0.25,1.75);
	\draw[very thick, directed=.65] (2,1) to [out=180,in=0] (-2,1);
	\draw[very thick, directed=.55] (2.5,0) to [out=180,in=0] (-1.5,0);
	\draw (2,1) to (2,5);
	\draw (2.5,0) to (2.5,4);
	\draw (-1.5,0) to (-1.5,4);
	\draw (-2,1) to (-2,5);	
	\draw[dashed] (2,3) to [out=180,in=45] (1,2.5);
	\draw[dashed] (2.5,2) to [out=180,in=315] (1,2.5);
	\draw[dashed] (1,2.5) to (-.5,2.5);
	\draw[dashed] (-.5,2.5) to [out=225,in=0] (-1.5,2);
	\draw[dashed] (-.5,2.5) to [out=135,in=0] (-2,3);
	\draw[very thick, red, directed=.65] (1,2.5) to [out=270,in=0]  (0.25,1.75) to [out=180, in = 270] (-0.5,2.5) to [out=90, in = 180] (0.25,3.25) to [out=0, in = 90] (1,2.5);
	\draw[very thick, directed=.65] (2,5) to [out=180,in=0] (-2,5);
	\draw[very thick, directed=.55] (2.5,4) to [out=180,in=0] (-1.5,4);
\end{tikzpicture}
\quad=\quad
\begin{tikzpicture} [anchorbase,scale=.5,fill opacity=0.2]
\path[fill=red]  (2,5) to [out=90,in=270] (2,1) to (-2,1) to (-2,5) to (2,5);
\path[fill=red]  (2.5,4) to [out=90,in=270] (2.5,0)  to (-1.5,0) to (-1.5,4) to (2.5,4);
	\draw[very thick, directed=.65] (2,1) to [out=180,in=0] (-2,1);
	\draw[very thick, directed=.55] (2.5,0) to [out=180,in=0] (-1.5,0);
	\draw (2,1) to (2,5);
	\draw (2.5,0) to (2.5,4);
	\draw (-1.5,0) to (-1.5,4);
	\draw (-2,1) to (-2,5);	
	\draw[dashed] (2,3) to (-2,3);
	\draw[dashed] (2.5,2) to (-1.5,2);
	\draw[very thick, directed=.65] (2,5) to [out=180,in=0] (-2,5);
	\draw[very thick, directed=.55] (2.5,4) to [out=180,in=0] (-1.5,4);
	\node[opacity=1] at (1.75,.5) {$\bullet$};
\end{tikzpicture}
\quad-\quad
\begin{tikzpicture} [anchorbase,scale=.5,fill opacity=0.2]
\path[fill=red]  (2,5) to [out=90,in=270] (2,1) to (-2,1) to (-2,5) to (2,5);
\path[fill=red]  (2.5,4) to [out=90,in=270] (2.5,0)  to (-1.5,0) to (-1.5,4) to (2.5,4);
	\draw[very thick, directed=.65] (2,1) to [out=180,in=0] (-2,1);
	\draw[very thick, directed=.55] (2.5,0) to [out=180,in=0] (-1.5,0);
	\draw (2,1) to (2,5);
	\draw (2.5,0) to (2.5,4);
	\draw (-1.5,0) to (-1.5,4);
	\draw (-2,1) to (-2,5);	
	\draw[dashed] (2,3) to (-2,3);
	\draw[dashed] (2.5,2) to (-1.5,2);
	\draw[very thick, directed=.65] (2,5) to [out=180,in=0] (-2,5);
	\draw[very thick, directed=.55] (2.5,4) to [out=180,in=0] (-1.5,4);
	\node[opacity=1] at (-1,4.5) {$\bullet$};
\end{tikzpicture}\]

In Khovanov homology as defined via Bar-Natan cobordisms, however, the base
point action changes by a sign when passing under a crossing in the link
diagram. To see this, consider the crossing complex
\[
  \CC{
    \begin{tikzpicture}[anchorbase,scale=.45]
      \draw [very thick, ->] (1,0) \pu (0,1.5);
      \draw [white, line width=.15cm] (0,0) \pu (1,1.5);
      \draw [very thick, ->] (0,0) \pu (1,1.5);
    \end{tikzpicture}
    }
    :=\left( 0 \to 
    \uwave{q^{1} \CC{
        \begin{tikzpicture}[anchorbase,scale=.45]
          \draw [very thick] (0,0) -- (0,1.5);
          \draw [very thick] (1,0) -- (1,1.5);
        \end{tikzpicture}
      }} \to q^{2}
      \CC{
        \begin{tikzpicture}[anchorbase,scale=.45]
          \draw [very thick] (0,0) \ur (.5,0.5) \rd (1,0);
          \draw [very thick] (0,1.5) \dr (.5,1) \ru (1,1.5);
        \end{tikzpicture} 
      } \to 0
      \right)
    \]
with differential given by the saddle cobordism. Again, there is a natural
homotopy given by the reverse saddle, but the commutator with the differential
is now the sum of the dots on both sides of the crossing, rather than the
difference. This is because of the famous neck cutting relation:
\begin{equation}
\label{eq:neckcut}
\begin{tikzpicture} [anchorbase,scale=.5,fill opacity=0.2]
\path[fill=red]  (2,5) to [out=90,in=270] (2,1) to [out=180,in=0] (0.25,1) to  (0.25,1.75) to [out=0, in=270] (1,2.5) to [out=90, in=0] (0.25,3.25) to [out=180, in=90] (-0.5,2.5) to [out=270,in=180] (0.25,1.75) to (0.25,1)  to (-2,1) to (-2,5) to (2,5);
\path[fill=red]  (2.5,4) to [out=90,in=270] (2.5,0) to [out=180,in=0] (0.25,0) to  (0.25,1.75) to [out=0, in=270] (1,2.5) to [out=90, in=0] (0.25,3.25) to [out=180, in=90] (-0.5,2.5) to [out=270,in=180] (0.25,1.75) to (0.25,0)  to (-1.5,0) to (-1.5,4) to (2.5,4);
	\draw[very thick] (2,1) to [out=180,in=0] (-2,1);
	\draw[very thick] (2.5,0) to [out=180,in=0] (-1.5,0);
	\draw (2,1) to (2,5);
	\draw (2.5,0) to (2.5,4);
	\draw (-1.5,0) to (-1.5,4);
	\draw (-2,1) to (-2,5);	
	\draw[dashed] (2,3) to [out=180,in=90] (1,2.5);
	\draw[dashed] (2.5,2) to [out=180,in=270] (1,2.5);
	\draw[dashed] (-.5,2.5) to [out=270,in=0] (-1.5,2);
	\draw[dashed] (-.5,2.5) to [out=90,in=0] (-2,3);
	\draw[very thick] (2,5) to [out=180,in=0] (-2,5);
	\draw[very thick] (2.5,4) to [out=180,in=0] (-1.5,4);
\end{tikzpicture}
\quad=\quad
\begin{tikzpicture} [anchorbase,scale=.5,fill opacity=0.2]
\path[fill=red]  (2,5) to [out=90,in=270] (2,1) to (-2,1) to (-2,5) to (2,5);
\path[fill=red]  (2.5,4) to [out=90,in=270] (2.5,0)  to (-1.5,0) to (-1.5,4) to (2.5,4);
	\draw[very thick] (2,1) to [out=180,in=0] (-2,1);
	\draw[very thick] (2.5,0) to [out=180,in=0] (-1.5,0);
	\draw (2,1) to (2,5);
	\draw (2.5,0) to (2.5,4);
	\draw (-1.5,0) to (-1.5,4);
	\draw (-2,1) to (-2,5);	
	\draw[dashed] (2,3) to (-2,3);
	\draw[dashed] (2.5,2) to (-1.5,2);
	\draw[very thick] (2,5) to [out=180,in=0] (-2,5);
	\draw[very thick] (2.5,4) to [out=180,in=0] (-1.5,4);
	\node[opacity=1] at (1.75,.5) {$\bullet$};
\end{tikzpicture}
\quad+\quad
\begin{tikzpicture} [anchorbase,scale=.5,fill opacity=0.2]
\path[fill=red]  (2,5) to [out=90,in=270] (2,1) to (-2,1) to (-2,5) to (2,5);
\path[fill=red]  (2.5,4) to [out=90,in=270] (2.5,0)  to (-1.5,0) to (-1.5,4) to (2.5,4);
	\draw[very thick] (2,1) to [out=180,in=0] (-2,1);
	\draw[very thick] (2.5,0) to [out=180,in=0] (-1.5,0);
	\draw (2,1) to (2,5);
	\draw (2.5,0) to (2.5,4);
	\draw (-1.5,0) to (-1.5,4);
	\draw (-2,1) to (-2,5);	
	\draw[dashed] (2,3) to (-2,3);
	\draw[dashed] (2.5,2) to (-1.5,2);
	\draw[very thick] (2,5) to [out=180,in=0] (-2,5);
	\draw[very thick] (2.5,4) to [out=180,in=0] (-1.5,4);
	\node[opacity=1] at (-1,4.5) {$\bullet$};
\end{tikzpicture}
\end{equation}

\subsection{Skein modules graded by integral homology}
The Kauffman bracket skein module $S_{\sltw}(M)$ of an oriented 3-manifold $M$
is defined as the quotient of the free module (over some chosen base ring)
spanned by all framed, unoriented links in $M$, modulo the Kauffman bracket
skein relations, i.e. the local relations satisfied by the Jones polynomial, see
Przytycki and Turaev \cite{MR1194712, MR1142906}. Links in $M$ represent classes
in $H_1(M,\Z/2\Z)$, endowing $S_{\sltw}(M)$ with a $H_1(M,\Z/2\Z)$-grading since
the skein relations are homogeneous. 

More generally, one can define skein modules $S_{\slN}(M)$,
which are naturally $H_1(M,\Z/N\Z)$-graded, or extend to the general linear case
and obtain skein modules $S_{\glN}(M)$ which are all $H_1(M,\Z)$-graded.

Kauffman bracket skein modules of thickened surfaces can be categorified using
Khovanov homology, see Asaeda--Przytycki--Sikora \cite{MR2113902}. A challenge
that arises in this context is that ordinary Khovanov homology (based on
unoriented Bar-Natan cobordisms) a priori does not prescribe chain maps
associated to non-orientable link cobordisms. However, such cobordisms are
unavoidable when working with ambient 3-manifolds that are general thickened
surfaces, and one common solution is to assign to non-orientable cobordisms the
zero map. This, however, has undesired consequences, since non-orientable
cobordisms are not closed under Bar-Natan's relations.

\newcommand{\torusfront}[3]{
\draw [green, thick, directed=.25] (0,0)to(#1,0) ;
\draw [red, thick, rdirected=.3, rdirected=.25] (#1+0.5*#2,#2) to (#1,0);
\draw [green, thick, directed=.25] (0,0+#3) to (#1,0+#3);
\draw [green, thick, directed=.25] (0+0.5*#2,#2+#3) to (#1+0.5*#2,#2+#3);
\draw [red, thick, rdirected=.3, rdirected=.25] (#1+0.5*#2,#2+#3) to (#1,0+#3);
\draw [red, thick, rdirected=.3, rdirected=.25] (0+0.5*#2,#2+#3) to (0,0+#3);
\draw (0,0) to (0,0+#3);
\draw(#1,0) to(#1,0+#3);
}

\newcommand{\torusback}[3]{
\draw [opacity=.5] (0.5*#2,#2)  to  (0.5*#2,#2+#3); 
\draw (#1+0.5*#2,#2) to (#1+0.5*#2,#2+#3);
\draw [green,opacity=.5, thick, directed=.25] (0+0.5*#2,#2)to (#1+0.5*#2,#2) ;
\draw [red,opacity=.5, thick, rdirected=.3, rdirected=.25] (0+0.5*#2,#2) to (0,0);
}
For example, in the thickened torus, consider the following unorientable cobordism between simple closed curves of distinct slope:

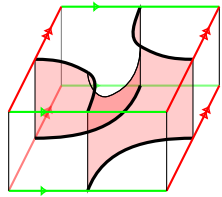
\[
\begin{tikzpicture}[fill opacity=.2,anchorbase,xscale=.7, yscale=0.7]
  \torusback{3}{2}{1.5}
  \draw [fill=red] (2.5,2.75) to [out=330,in=180] (3.5,2.5) to (3.5,1)to [out=180,in=70](1.5,0) to (1.5,1.5) to [out=70,in=330] (1.5,2.25) to [out=270,in=200] (1.9,1.8) to [out=20,in=270] (2.5,2.75);
  \draw [fill=red] (2.5,2.75) to [out=150,in=250](2.5,3.5) to (2.5,2) to [out=250,in=0] (.5,1) to  (.5,2.5) to [out=0,in=150] (1.5,2.25) to [out=270,in=200] (1.9,1.8) to [out=20,in=270] (2.5,2.75) ; 
  \draw[very thick] (.5,1) to [out=0,in=250](2.5,2);
  \draw[very thick] (1.5,0) to [out=70,in=180](3.5,1);
  \draw[very thick] (.5,2.5) to [out=0,in=150] (1.5,2.25);
  \draw[very thick] (1.5,1.5) to [out=70,in=330] (1.5,2.25);
  \draw[very thick] (2.5,2.75) to [out=330,in=180](3.5,2.5);
  \draw[very thick] (2.5,2.75) to [out=150,in=250](2.5,3.5);
  \torusfront{3}{2}{1.5}
  \end{tikzpicture}
  \]
  The composite of this cobordism with its reverse can be neck-cut as in
  \eqref{eq:neckcut}, to obtained a dotted (orientable) identity cobordism
  scaled by $2$. Setting unorientable cobordisms to zero, therefore, also kills
  certain orientable cobordisms. For example, in \cite{MR2475122} this is
  accounted for by the relation (NDD).

Foams circumvent these phenomena by requiring all facets to be oriented.
Specifically, when working with $\gltw$ foams, one can show that the class of
foams with non-orientable underlying 1-labelled surfaces are closed under local
Blanchet foam relations, see Queffelec--Wedrich~\cite[Proposition
3.11]{1806.03416}. This leads to a categorification of the $\gltw$ skein modules
of thickened surfaces in \cite{1806.03416}, where we also proposed candidates
for categorified $\glN$ skein modules using Khovanov--Rozansky homology. 

One dimension higher, Khovanov--Rozansky homology can be used to build skein
modules $\skeinzero(W^4;L)$ based on framed, oriented surfaces embedded in
smooth oriented 4-manifolds $W^4$, possible with a link $L$ in the boundary, see
Morrison--Walker--Wedrich~\cite{2019arXiv190712194M}. The compatible orientation
of facets in $\glN$ foams and the homogeneity of the foam relations implies that
these skein modules are naturally graded by $H_2(W^4;L,\Z)$ in addition to the
cohomological and quantum grading. A version based on $\slN$ would be expected to carry a grading by $H_2(W^4;L,\Z/N\Z)$, see e.g. the discussion acound \cite[Equation (14)]{ren2024khovanov} in the case of $N=2$. 

Skein modules based on $\glN$ link homology, furthermore, admit deformations which decompose along $\glN$ branching rules, generalizing \eqref{eq:coldef}, see \cite{morrison2024invariantssurfacessmooth4manifolds}.

\renewcommand*{\bibfont}{\small}
\setlength{\bibitemsep}{0pt}
\raggedright
\printbibliography

\end{document}